\documentclass[review,
]{elsarticle}

\usepackage{graphicx,color}

\usepackage[colorlinks=true, breaklinks=true, pdfstartview=FitV, linkcolor=red, citecolor=blue, urlcolor=black]
{hyperref}

\usepackage[T1]{fontenc}

\journal{Acta Astronautica; doi: 10.1016/j.actaastro.2019.01.039\quad}

\begin{document}

\begin{frontmatter}

\title{ Design of distant retrograde orbits based on a higher order analytical solution\tnoteref{t1} }
\tnotetext[t1]{Presented at the 69th International Astronautical Congress 2018, Bremen, Germany, 1-5 October 2018, as paper IAC-18.C1.1.1.x42642}

\author[RA]{Martin Lara
\fnref{footnote2,footnote1}}
\ead{mlara0@gmail.com}

\address[RA]{Edificio CCT, C/ Madre de Dios, 53, ES-26006 Logro\~no, Spain}

\fntext[footnote2]{GRUCACI, University of La Rioja}
\fntext[footnote1]{Space Dynamics Group, Technical University of Madrid -- UPM}


\date{}
\begin{abstract}
Because the distant retrograde orbits dynamics inherently depends on special functions, approximate analytical solutions in the literature are commonly constrained to providing rough approximations of the qualitative behavior. We rely on perturbation methods and succeed in computing a higher order analytical solution that extends the range of applicability of previous solutions to the problem of relative, quasi-satellite orbits with large librations. Besides, the analytical solution provides two design parameters that are effectively used in the computation of periodic, distant retrograde orbits without constraining to the typical 1:1 resonance.
\end{abstract}


\end{frontmatter}

\section{Introduction}

New proposed space missions, as the Deimos and Phobos Interior Explorer proposal to ESA's Cosmic Vision Program \cite{OberstWillnerWickhusen2017}, or NASA's Asteroid Redirect Mission \cite{Abelletal2017}, motivate current interest in the study of distant retrograde orbits (DRO). 
\par

Proposed solutions to the DRO problem usually rely on the numerical computation of periodic and quasi-periodic orbits of a restricted three-body problem model \cite{Benest1974,LidovVashkovyak1994,LamWhiffen2005,LaraRussellVillac2007,MingShijie2009,GilSchwartz2010,StramacchiaColomboBernelliZazzera2016,BezroukParker2017}. Analytical approximations to the DRO problem can also be computed, but their accuracy is commonly constrained to rough approximations of the dynamics, which are useful only in the qualitative descriptions of the co-orbital dynamics \cite{Benest1976,Namouni1999}.
\par

In spite of the performance of analytical DRO solutions can be improved when approaching the problem from the point of view of Hamiltonian perturbation methods \cite{LidovVashkovyak1993}, limitations in the accuracy of perturbation solutions due to the inherent dependence of the dynamics on special functions, which limit the practical implementation of the perturbation solution to lower orders, made researchers to consider the perturbation approach just partially effective \cite{LidovVashkovyak1994a}. In consequence, the use of these kinds of solutions was limited to the description of the long term orbit behavior, which was obtained by numerically integrating the evolution equations that are obtained after the elimination of short period effects \cite{LidovVashkovyak1994}.
\par

It was, perhaps, the untimely demise of Lidov which made the perturbation approach to be abandoned until recently. Still, new applications of the perturbation approach are, again, limited to the orbit evolution, which equations are obtained by numerically averaging the higher frequencies of the motion \cite{Sidorenkoetal2014}. However, purely analytical efforts can be done without the need of relying on numerical averaging. In particular, focusing on the approximation to the dynamics provided by the Hill problem, both a low and a high order perturbation solutions have been recently reported without limitation to the evolutionary equation \cite{Lara2018}. Indeed, further than the usual averaging, the short-period corrections that provide the transformation from osculating to mean elements and vice-versa, have also been obtained analytically, yet only to some order of the perturbation approach due to the special functions that arise in the procedure when the coupling of the different disturbing effects takes place in the perturbation approach.
\par

Thus, after formulating the relative motion described by the Hill problem in epicyclic variables, the phase of the orbiter in its ``epicycle'' is removed by perturbation methods in \cite{Lara2018}. As a result of the elimination of this fast angle, and up to certain truncation order of the perturbation approach, the planar Hill problem is decoupled into a reduced system of one degree of freedom for the motion of the ``deferent'', which is integrable, and a quadrature that provides the epicycle.\footnote{At difference from the planetary motion, both the epicycle and the deferent are ellipses, contrary to circumferences, in the case of DRO, and the size of the deferent is smaller than that of the epicycle.}
\par

For the lower orders of the DRO perturbation approach, the flow can be integrated with an exact, explicit, analytical solution. This solution only depends on trigonometric functions and, therefore, has the advantage of bringing deeper insight on the dynamics, as well as providing faster evaluation, than previous alternatives in the literature ---which depend on special, contrary to trigonometric functions. Notably, the low order solution obtained with the perturbation approach discloses the existence of two basic design parameters that can be used for the computation of initial conditions of orbits that are periodic on average, and almost periodic in the original Hill problem. Following application of iterative differential corrections to the initial conditions predicted by the analytical solution makes them to easily converge to the initial conditions and period of a true periodic orbit with the characteristics fixed by the design parameters ---a procedure that is standard in preliminary mission design \cite{SADSaM1999,Lara1999,SanJuanLaraFerrer2006,Lara2011}.
\par

The computation of higher orders of the perturbation approach extends the range of applicability of the analytical solution in epicyclic variables to the case of DRO with larger librations. However, the advantages of getting an exact solution in trigonometric functions are inevitably lost. Indeed, while the reduced flow is still integrable the solution is achieved in an implicit way, on the one hand, and rely on the use of special functions, on the other \cite{Lara2018}. While both inconveniences degrade the performance of the perturbation solution, they can be partially remedied by the use of Lindstedt series. This procedure makes the computation of an explicit approximation to the higher order solution that still relies only on trigonometric functions again feasible.
\par

The Lindstedt series solution computed in \cite{Lara2018} improved applicability of the analytical solution to the DRO problem with respect to other similar approaches in the literature. However, it was constrained to the first terms of both the perturbation approach and the Lindstedt series. Therefore, it failed in the accurate modeling of DRO with large librations, for which case the usual semi-analytical integration of the averaged equations, which progresses with very large step sizes, was suggested as the way to obtain more accurate solutions in a fast an efficient way \cite{LidovVashkovyak1994}. 
\par

Here we show how the computation of additional terms of the Lindstedt series solution overcomes some of the limitations of previous analytical approaches, in this way allowing us to model higher order effects of the perturbation due to the gravitation of the central body  ---the co-orbital companion of the massless body in its journey about the larger primary--- in a pure analytical way. The applicability of the analytical solution, which now can account for wider amplitudes of the orbit libration, is then notably broaden, and the accuracy of the higher order Lindstedt series solution is found to be comparable to that of the semi-analytical integration of DRO with large librations, with the advantage on the side of the Lindstedt series solution of avoiding the need of carrying out any numerical integration.
\par

The paper is organized as follows. In Section 2, we recall basic facts on the formulation of the Hill problem in epicyclic variables, which are traditionally used in the study of quasi-satellite orbits. While this formulation has been previously used as well as thoroughly described by different authors, we depart from the tradition by introducing a scaling parameter in the standard transformation to epicyclic variables. Because of that, we feel compelled to give some detail into the somewhat different equations to be solved. Section 3 deals with the short-period averaging, which yields the (truncated) long-term Hamiltonian from which the evolution equations stem. Due to the scaling made in the canonical transformation to epicyclic variables, the long-term Hamiltonian of the DRO problem is clearly cast in the form of two coupled perturbed harmonic oscillators. Details on the nature of the perturbed solution are provided in Section 4, whereas the construction of an improved solution using Lindstedt series is discussed in Section 5. Finally, Section 6 provides different examples that show the performance of the new analytical solutions, and illustrate the use of the two design parameters in the construction of single- and multiple-period periodic orbits.
\par

\section{The planar Hill problem as a perturbed quadratic Hamiltonian}

The relative motion of two bodies that move in co-orbital motion about a distant, massive primary can be efficiently modeled in the Hill problem approximation \cite{Hill1878}. Thus, in the restricted approximation in which one of the bodies has negligible mass, the primaries are assumed to evolve with Keplerian motion, while the orbit of the body of negligible mass (the orbiter) is attracted by both primaries. If the motion of the primaries is further assumed to be circular, and the distance of the orbiter to the lesser mass primary is always small compared to the distance between both primaries, the relative motion of the orbiter with respect to the less massive primary can be derived from the Hamiltonian
\begin{equation} \label{planarHam}
\mathcal{H}=\mathcal{Q}(x,y,X,Y;\omega)-\frac{\mu}r,
\end{equation}
in which $\mu$ is the gravitational parameter of the lesser mass primary, $r=\sqrt{x^2+y^2}$ and
\begin{equation} \label{H0}
\mathcal{Q}=\frac{1}2\left(X+\omega{y}\right)^2+\frac{1}2(Y-\omega{x})^2-\frac{3}{2}\omega^2x^2,
\end{equation}
is a quadratic Hamiltonian.
\par

Note that we constrain to the case of planar motion, in the orbital plane of the primaries, and the Cartesian coordinates of the orbiter $(x,y)$ as well as their conjugate momenta $(X,Y)$ are referred to a direct rotating frame with rotation rate $\omega$, with the $x$ axis direction defined by the line joining both primaries in the direction from the more massive to the less massive one. Full details on the formulation of the problem can be consulted in textbooks on celestial mechanics, as, for instance, \cite{Szebehely1967}.
\par

The physical parameters $\mu$ and $\omega$ upon which the Hill problem seems to depend, are in fact superfluous because they can be ignored by the usual choice of Hill units in which $\omega=\mu=1$. The absence of physical parameters gives a broad generality to the Hill problem. However, keeping both quantities, $\mu$ and $\omega$, explicit is quite useful in making analytical derivations for the insight they provide into the problem, but also because they permit the systematic check of dimensions as an additional test to be applied at each step of the procedure. So we hold formally both parameters in what follows.
\par

If the orbiter is either far away from the central body or the gravitational parameter of this primary is very small, then the contribution of the summand $\mu/r$ in Eq.~(\ref{planarHam}) can be neglected. Then, the Hamiltonian defined by the remaining terms, that is, by Eq.~(\ref{H0}), is quadratic, and, therefore, the associated Hamilton equations form a linear differential system ---the so called Clohessy-Wiltshire equations \cite{ClohessyWiltshire1960}--- whose solution is known to be a drifting ellipse in the $y$ axis direction with a ratio semi-minor semi-major axis 1:2, for the orbit, whereas the hodograph is a drifting circle in the $X$ axis direction. This drift can be nullified by the adequate choice of initial conditions, which, therefore provide periodic solutions.
\par

Instead of using the classical solution, and in preparation for a Hamiltonian perturbation approach, we better find the solution to the flow stemming from Eq.~(\ref{H0}) by complete reduction, which is achieved by the canonical transformation to epicyclic variables $T:(x,y,X,Y)\rightarrow(\phi,q,\Phi,Q;\omega,k)$, given by
\begin{equation} \label{tdir}
\begin{array}{rcl}
x &=& 2b\xi+b\sin\phi,  \\
y &=& a\eta+a\cos\phi, \\
X &=& -2B\eta-B\cos\phi, \\
Y &=& -B\xi-B\sin\phi,
\end{array}
\end{equation}
with $a=2b$, $B=b\omega$, and
\begin{equation} \label{bxy}
b=\sqrt{2\Phi/\omega}, \quad \xi=\frac{Q}{2kB}, \quad \eta=\frac{2kq}{a}.
\end{equation}
The value of the scaling factor $k$ is yet to be assigned. 
\par

Equation (\ref{tdir}) can be seen as the parametric representation of an ellipse with semi-axes $a$ and $b$, and center
\begin{equation} \label{center}
x_\mathrm{C}=2b\xi=Q/(k\omega), \qquad y_\mathrm{C}=a\eta=2kq.
\end{equation}
for the orbit, and a circumference of radius $B$ and center $(-2B\eta,-B\xi)$ for the hodograph.
\par

The transformation given by Eqs.~(\ref{tdir})--(\ref{bxy}) turns Eq.~(\ref{H0}) into a Hamiltonian which is completely reduced in the new variables, viz.
\begin{equation} \label{quadratic}
\mathcal{Q}=\mathcal{Q}(-,-,\Phi,Q)\equiv\omega\Phi-(3/8)(Q/k)^2,
\end{equation}
in which $\Phi$ and $Q$ are constant because their conjugate coordinates $\phi$ and $q$, respectively, are cyclic. The constant values of $\Phi$ and $Q$ depend on the particular choice of the initial conditions $(x_0,y_0,X_0,Y_0)$, that is
\begin{eqnarray}
\Phi &=& \frac{1}{2\omega}\left[(X_0+\omega{y}_0)^2+(2Y_0+\omega{x}_0)^2\right], \\
Q &=& 2k(Y_0+\omega{x}_0).
\end{eqnarray}
\par

The Hamiltonian flow originated from Eq.~(\ref{quadratic}) enjoys trivial solution. Indeed,
\begin{equation} \label{fiqp}
\dot\phi=\frac{\partial\mathcal{Q}}{\partial\Phi}= \omega, \qquad
\dot{q}=\frac{\partial\mathcal{Q}}{\partial{Q}}=-\frac{3/4}{k^2}Q,
\end{equation}
where $\omega$ and $-\frac{3}{4}Q/k^2$ are constant frequencies. The latter suggests to choose
\begin{equation} \label{kchoice}
k=\sqrt{3/4},
\end{equation}
This choice of the scaling parameter is applied throughout the remaining of the paper yet we find convenient to keep the notation $k$ in some places in order to avoid as much as possible the explicit appearance of square roots.\footnote{Note that this choice of $k$ is different from the value $k=1$ used in \cite{Lara2018}, with the consequent changes in derived expressions. This is the reason why we note the auxiliary variables differently from \cite{Lara2018}, in this way trying to avoid confusion in prospective readers of both works.}
Hence
\begin{eqnarray} \label{solphi}
\phi &=& \phi(0)+\omega{t}, \\  \label{solq} 
q &=& q(0)-Qt.
\end{eqnarray}
\par

Therefore, when $Q=0$ the orbit is just an ellipse. On the contrary, both $y$ and $X$ will grow unbounded in a linear way when $Q\ne0$, but still the orbit can be viewed as a moving ellipse whose guiding center in Eq.~(\ref{center}) evolves with linear motion in the $y$ axis direction.
\par

If we also apply Eq.~(\ref{tdir}) to convert the term $-\mu/r$, which makes the Hamiltonian (\ref{planarHam}) nonlinear, into epicyclic coordinates we get
\begin{equation} \label{muir}
\frac{\mu}{r}=\frac{\mu}{a}\frac{1}{\sqrt{\Delta^2+\xi\sin\phi+2\eta\cos\phi+\xi^2+\eta^2}},
\end{equation}
with 
\begin{equation} \label{sigma}
\Delta=\sqrt{1-k^2\sin^2\phi}.
\end{equation}
\par

For conciseness, in what follows we use the abbreviations $c\equiv\cos\phi$, $s\equiv\sin\phi$. Then, the Hill problem Hamiltonian (\ref{planarHam}) is written in epicyclic coordinates as
\begin{equation} \label{Hameq}
\mathcal{H}=\omega\Phi\bigg(1-3\xi^2-\frac{\gamma}{\sqrt{\Delta^2+\xi{s}+2\eta{c}+\xi^2+\eta^2}}\bigg),
\end{equation}
in which $\gamma\equiv\gamma(\Phi)$ is the non-dimensional, auxiliary variable
\begin{equation} \label{gamma}
\gamma=\frac{\mu}{a\omega\Phi}=\frac{\mu\omega}{(2\omega\Phi)^{3/2}}.
\end{equation}

\section{Long-term Hamiltonian}

Solution of the DRO problem by perturbations is eased when expanding the inverse of the square root in Eq.~(\ref{Hameq}) in power series of some small quantity. To make things more clear, we introduce a formal small parameter $\epsilon$ and scale the problem making the assumptions
\begin{equation} \label{assumptions}
\eta=\mathcal{O}(\epsilon), \qquad
\xi=\mathcal{O}(\epsilon^2), \qquad
\gamma=\mathcal{O}(\epsilon^4).
\end{equation}
That is, we implicitly assume that the 3rd body and central gravitation perturbations are of the same order, cf.~\cite{Lara2018}. According to these assumptions, we arrange the expanded Hamiltonian in the form
\begin{equation} \label{Hamp}
\mathcal{H}=\sum_{j\ge0}\frac{\epsilon^j}{j!}H_j,
\end{equation}
The first few terms of Eq.~(\ref{Hamp}) are $H_0=\omega\Phi$, $H_1=H_2=H_3=0$, and
\begin{eqnarray*}
\frac{H_4}{H_0} &=&\!\!\!\! 4!\,(-1)\Big(3\xi^2+\frac{\gamma}{\Delta}\Big)
\\
\frac{H_5}{H_0} &=& 5!\,\frac{\gamma}{\Delta^3}\eta{c}, \\
\frac{H_6}{H_0} &=& 6!\,\frac{\gamma}{8 \Delta ^5}\big[\eta^2(1-9c^2)+\xi\left(1+3c^2\right)s\big] \\
\frac{H_7}{H_0} &=& 7!\,\frac{\gamma}{8\Delta^7}\big[\eta^2(8-11s^2)-3\xi(1+3c^2)s\big]\eta{c} \\
\frac{H_8}{H_0} &=& 8!\,\frac{\gamma}{\Delta^9}\Big[
\Big(-1+\frac{11}{4}s^2-\frac{227}{128}s^4\Big)\eta^4 \\
&&+\eta^2\xi\Big(3-\frac{87}{16}s^2+\frac{153}{64}s^4\Big)s \\
&&+\xi^2\Big(\frac{1}{2}-\frac{3}{2}s^2+\frac{45}{32}s^4-\frac{27}{64}s^6\Big)
\Big]
\end{eqnarray*}
\par

The perturbation Hamiltonian in Eq.~(\ref{Hamp}) is truncated to the order 13th, because of the known difficulties in progressing to higher orders in the elimination of the fast angle $\phi$ introduced by the gradual coupling of the different orders of the perturbation approach, cf.~\cite{Lara2018}.
\par

After removing $\phi$ by a canonical transformation to new, prime variables $(\phi,q,\Phi,Q)\rightarrow(\phi',q',\Phi',Q';\epsilon)$, we obtain
\begin{equation} \label{H8a}
\mathcal{H}=\omega\Phi'-\frac{1}{2}(Q'^2+\Omega^2q'^2)+\mathcal{P}(-,q',\Phi',Q'),\
\end{equation}
where the \emph{libration frequency}
\begin{equation} \label{lowfreq}
\Omega=\Omega(\Phi')\equiv\omega\sqrt{\tilde\mathrm{K}-\tilde\mathrm{E}}\sqrt{\gamma},
\end{equation}
in which $\gamma$ is given by Eq.~(\ref{gamma}) with $\Phi$  replaced by $\Phi'$, and we abbreviated $\tilde\mathrm{K}\equiv\mathrm{K}(k^2)/\pi$, $\tilde\mathrm{E}\equiv\mathrm{E}(k^2)/\pi$, with $\mathrm{K}$ and $\mathrm{E}$ denoting the complete elliptic integrals of the 1st and 2nd kind, respectively. The term $\mathcal{P}$ is
\begin{equation} \label{perturbation}
\mathcal{P}=\omega\Phi'\sum_{i=0}^2\sum_{j=0}^{2-i}\sum_{n=0}^{4-2i-2j}p_{i,j,n}\gamma^{i+1}\xi^{2j},\eta^{2n},
\end{equation}
which numeric coefficients $p_{i,j,n}$ are given in Table \ref{t:pijn}.
\par

Note that $\xi$ and $\eta$ in Eq.~(\ref{perturbation}) are still given by Eq.~(\ref{bxy}), but now $q$, $Q$, and $\Phi$, must be replaced by $q'$, $Q'$, and $\Phi'$, respectively. Besides since we assumed $\gamma=\mathcal{O}(\epsilon^4)$ in Eq.~(\ref{assumptions}), it happens that $\Omega/\omega=\mathcal{O}(\epsilon^2)$, as follows from Eq.~(\ref{lowfreq}). Finally, it is worth noting that, due to the Hamiltonian arrangement, the coupling between the different orders is avoided up to the 8th order of the perturbation approach. Hence, the first 7 new Hamiltonian terms are just the average of the old terms over $\phi$.
\par

\begin{table}[htb]
\caption{Coefficients $p_{i,j,n}$ in Eq. (\protect\ref{perturbation})}
\centering
\begin{tabular}{ll}\hline
$p_{0,0,0}=-2\tilde\mathrm{K}$ \vphantom{$\frac{M^2}{0}$} \\[0.5ex]
$p_{0,0,1}=0$ \\
$p_{0,0,2}=\frac{1}{9}(14\tilde\mathrm{E}-11\tilde\mathrm{K})$ \\[0.5ex]
$p_{0,0,3}=\frac{2}{81}(71\tilde\mathrm{E}-50\tilde\mathrm{K})$ \\[0.5ex]
$p_{0,0,4}=\frac{1}{324}(644\tilde\mathrm{E}-425\tilde\mathrm{K})$ \\[0.5ex]
$p_{0,1,0}=\frac{4}{3}(\tilde\mathrm{K}-4\tilde\mathrm{E})$ \\[0.5ex]
$p_{0,1,1}=4(\tilde\mathrm{K}-12\tilde\mathrm{E})$ \\[0.5ex]
$p_{0,1,2}=-\frac{10}{27}(74 \tilde\mathrm{E}-35 \tilde\mathrm{K})$ \\[0.5ex]
$p_{0,2,0}=\frac{4}{9}(\tilde\mathrm{K}-16\tilde\mathrm{E})$ \\[0.5ex]
$p_{1,0,0}=\frac{1}{2}-2\tilde\mathrm{K}^2$ \\[0.5ex]
$p_{1,0,1}=\frac{16}{9}\big(\frac{21}{32}+2\tilde\mathrm{E}^2+\tilde\mathrm{E}\tilde\mathrm{K}-3\tilde\mathrm{K}^2\big)$ \\[0.5ex]
$p_{1,1,0}=\frac{8}{3}(3-8\tilde\mathrm{E}\tilde\mathrm{K}+2\tilde\mathrm{K}^2)$ \\[0.5ex]
$p_{1,0,2}=\frac{1}{3}(9-8\tilde\mathrm{E}^2+44\tilde\mathrm{E}\tilde\mathrm{K}-30\tilde\mathrm{K}^2)$ \\[0.5ex]
$p_{2,0,0}=5\big(\frac{4}{15}-\tilde\mathrm{E}+\tilde\mathrm{K}-\tilde\mathrm{K}^3\big)$ \\[0.5ex]
\hline
\end{tabular}
\label{t:pijn}
\end{table}

Note that, on average, the motion of the center of the reference ellipse decouples from the $\phi'$ motion. Indeed, since $\phi'$ is ignorable in Eq.~(\ref{H8a}), $\Phi'$ is a dynamical parameter which remains constant. Therefore, the reduced $(q',Q')$ system can be integrated from Hamilton equations
\begin{equation} \label{reduced}
\frac{\mathrm{d}q'}{\mathrm{d}t}=\frac{\partial\mathcal{H}}{\partial{Q}'}, \qquad
\frac{\mathrm{d}Q'}{\mathrm{d}t}=-\frac{\partial\mathcal{H}}{\partial{q}'}.
\end{equation}
Once the solution of the reduced system (\ref{reduced}) is known, viz.~$q=q(t,q'_0,Q'_0;\Phi')$, $Q=Q(t,q'_0,Q'_0;\Phi')$, the remaining variable $\phi'$ is integrated by quadrature
\begin{equation} \label{quadrature}
\phi'=\phi'_0+\int\frac{\partial\mathcal{H}(t,q_0,Q_0;\Phi)}{\partial{\Phi}}\,\mathrm{d}t.
\end{equation}
\par

The solution of Eqs.~(\ref{reduced})--(\ref{quadrature}) just provides the DRO evolution in the prime, \emph{mean} variables. The \emph{osculating} motion is recovered after applying the canonical transformation from prime to original variables. Explicit expressions of this transformation are provided in \ref{a:peco}.

\section{The nature of the solution}

The main characteristics of DRO of the Hill problem are disclosed after neglecting $\mathcal{P}$ from Eq.~(\ref{H8a}), showing that, on average, the motion is composed of two perturbed harmonic oscillations. The Hamiltonian of the first (unperturbed) harmonic oscillator is $\mathcal{F}=\omega\Phi'$, and gives the phase of \emph{epicyclical} motion of the orbiter $\phi'=\phi'_0+\omega{t}$, with respect to the center of the reference ellipse. The second (unperturbed) harmonic oscillator is $\mathcal{G}=-\frac{1}{2}(Q'^2+\Omega^2q'^2)$, and gives the \emph{deferential} motion of the center of the reference ellipse with respect to the the primary of lesser mass. 
\par

An approximate explicit solution to Eq.~(\ref{H8a}) is easily computed if we neglect the contribution of $\mathcal{P}$. In fact, to be consistent with the assumptions (\ref{assumptions}) made for the arrangement of the perturbation Hamiltonian, we must keep from $\mathcal{P}$ the first term of the summation (\ref{perturbation}), viz. $\omega\Phi'\gamma{p}_{0,0,0}$. Then, we only deal with the Hamiltonian
\begin{equation} \label{H0p}
\mathcal{H}_0=\tilde\omega\Phi'-\frac{1}{2}(Q'^2+\Omega^2q'^2),
\end{equation}
where $\tilde\omega=\omega(1+\gamma{p}_{0,0,0})$, and we recall that both $\gamma$ and $\Omega$ are function of $\Phi'$.\footnote{An attentive reader may have noticed that including additional terms of the perturbation into Eq.~(\ref{H0p}) yields the Duffing oscillator, in this way getting a more comprehensive integrable zero order Hamiltonian. However, the solution of the Duffing oscillator needs the use of Jacobi elliptic functions, thus lacking of the insight that we pursue. }

\subsection{Deferential harmonic oscillations}

Indeed, to solve the reduced flow from Eq.~(\ref{reduced}) we apply the usual Poincar\'e transformation $(q',Q')\rightarrow(g,G;\Omega)$ defined by
\begin{equation} \label{poin}
q'=\sqrt{2G/\Omega}\sin{g}, \qquad Q'=\sqrt{2\Omega{G}}\cos{g},
\end{equation}
which transforms Eq.~(\ref{H0p}) into $\mathcal{H}_0=\tilde\omega\Phi'-\Omega{G}$. Then, in addition to $\Phi'$, $G$ is also constant and, therefore, $g=g_0-\Omega{t}$, where $g_0$ is unambiguously determined from the inverse transformation of Eq.~(\ref{poin}) as a function of $q'_0$ and $Q'_0$.
\par

The solution in the $(q',Q')$ variables is recovered from Eq.~(\ref{poin}), leading to
\begin{eqnarray} \label{qosc}
q' &=& q'_0\cos\Omega{t}-(Q'_0/\Omega)\sin\Omega{t}, \\ \label{qposc}
Q' &=& Q'_0\cos\Omega{t}+\Omega{q}'_0\sin\Omega{t},
\end{eqnarray}
which in turn is replaced into Eq.~(\ref{quadrature}) to yield
\begin{equation} \label{fi4th}
\phi'=\phi'_0+\omega\left(1+\frac{\Omega^2}{\omega^2}d\right)t+\frac{\Omega}{\omega}\big[p(t)-p(0)\big],
\end{equation}
with
\begin{equation} \label{dd0}
d=\frac{\tilde\mathrm{K}}{\tilde\mathrm{K}-\tilde\mathrm{E}}+\frac{{q'_0}{}^2+(Q'_0/\Omega)^2}{(b/k)^2}>0,
\end{equation}
and
\begin{equation} \label{pp0}
p=\frac{{q}'_0(Q'_0/\Omega)}{(b/k)^2}\cos2\Omega{t}+\frac{{q'_0}{}^2-(Q'_0/\Omega)^2}{2(b/k)^2}\sin2\Omega{t}.
\end{equation}
\par

Note that, except for the $k$ scaling in the transformation to epicyclic variables in Eqs.~(\ref{tdir})--(\ref{bxy}), the solution given by Eqs.~(\ref{qosc})--(\ref{pp0}) is the same as the solution in Eqs.~(43)--(46) of \cite{Lara2018}. That is, the center of the reference ellipse in Eq.~(\ref{center}) evolves slowly with harmonic oscillations
\begin{equation} \label{centeraverage}
x_C=\frac{\Omega}{k\omega}M\sin(\Omega{t}+\psi), \quad y_C=2kM\cos(\Omega{t}+\psi),
\end{equation}
with
\begin{equation} \label{Mpsi}
M=\sqrt{q'_0{}^2+(Q'_0/\Omega)^2}, \quad \tan\psi=Q'_0/(\Omega{q}'_0).
\end{equation}
Thus, on average, the center of the reference ellipse moves along another ellipse with the semi-major axes also in the $y$ axis direction and notably longer tha the semi-minor axis because of the factor $\Omega/\omega$ which is $\mathcal{O}(\epsilon^2)$ in our perturbation assumptions. On the other hand, the linear growing of the phase of the satellite, which happens at a rate slightly faster than the rotation rate of the system $\omega$, is perturbed by harmonic oscillations of small amplitude and twice the frequency with which the center of the reference ellipse evolves.
\par

\subsection{Orbit design parameters}
The low-order solution comprising Eqs.~(\ref{qosc})--(\ref{pp0}) shows that DRO are essentially defied by two basic parameters: the size of the drifting ellipse, as given by any of the semi-axes, and the minimum distance in the $y$ axis direction to the smaller primary $\rho=a-2kM$, cf.~\cite{LidovVashkovyak1994}. \par

Indeed, given $b$, $\Phi$ is computed from the first of Eq.~(\ref{bxy}). Then, $\Omega$ is obtained using Eqs.~(\ref{lowfreq}) and (\ref{gamma}). Initial conditions of the deferential motion are then solved from Eq.~(\ref{Mpsi}), in which $M=(a-\rho)/(2k)$ and $\psi$ is arbitrarily chosen. Thus, if, for instance, we chose $\psi=\pi/2$, then $q_0=0$ and $Q_0=\Omega(a-\rho)/(2k)$. If we choose $\psi=0$, then $Q_0=0$, and $q_0=(a-\rho)/(2k)$. 
\par

On the other hand, we are free to choose any value of $\phi$ to initialize the epicyclical motion, in this way completing an initial set of initial conditions of a DRO with the desired characteristics in the averaged, prime space.
\par

\section{Improving accuracy: Lindstedt's series approach}

When the remaining terms of the perturbation $\mathcal{P}$ are taken into account, the long-term Hamiltonian remains of one degree of freedom, but an explicit analytical solution to it has not been found yet. Alternatively, the use of approximate solutions based on truncated Lindstedt series has shown successful in the analytical computations of DRO with small to moderate librations. However, the accuracy of the proposed series notably degrades in the case of larger librations. Here we extend the Lindstedt series solution of \cite{Lara2018} to encompass higher order effects of the perturbation, and show that the new solution provides results of comparable accuracy to the semi-analytical integration that is customarily approached for DRO with large librations \cite{LidovVashkovyak1994,Lara2018}.
\par

The basics of the Lindstedt series approach is as follows. First, we change the independent variable
\begin{equation}
\tau=nt,
\end{equation}
in which the Lindstedt series are solved. The effect of this change is a simple time scaling of Eq.~(\ref{reduced}), which now reads
\begin{equation} \label{redtau}
n\frac{\mathrm{d}q'}{\mathrm{d}\tau}=\frac{\partial\mathcal{H}'}{\partial{Q}'}, \qquad
n\frac{\mathrm{d}Q'}{\mathrm{d}\tau}=-\frac{\partial\mathcal{H}'}{\partial{q}'}.
\end{equation}
Next, the series
\begin{equation}
\begin{array}{rcl}
n &=& \displaystyle 1+\sum_{i\ge1}\varepsilon^in_i, \\[2.5ex]
q' &=& \displaystyle \sum_{i\ge0}\varepsilon^iq^{(i)}(\tau), \\[2.5ex]
Q' &=& \displaystyle \sum_{i\ge0}\varepsilon^iQ^{(i)}(\tau),
\end{array}
\end{equation}
where $\varepsilon$ is just formal, are replaced into Eq.~(\ref{redtau}). Then, after equating the coefficients of equal powers of $\varepsilon$, we obtain a chain of differential systems that can be solved sequentially as a function of $\tau$ and the initial conditions $q'_0$, $Q'_0$. Finally, for each system $i$, we must determine the coefficient $n_i$ so that the solution $(q^{(i)}(\tau),Q^{(i)}(\tau))$ be periodic. Full details on the computation of Lindstedt series can be consulted in textbooks on perturbation methods like, for instance, \cite{Murdock1999,Nayfeh2004}.
\par

Once $q'$ and $Q'$ have been solved in the new time scale up to the desired order of $\varepsilon$, the long-term evolution of the phase of the satellite is computed from Eq.~(\ref{quadrature}), which in the new time scale reads
\begin{equation} \label{quadtau}
\phi'=\phi'_0+\frac{1}{n}\int\frac{\partial\mathcal{H}'(q'(\tau),Q'(\tau);\Phi')}{\partial{\Phi'}}\,\mathrm{d}\tau.
\end{equation}
This integral can also be solved sequentially by computing at each step of the procedure the corresponding term $\phi_i(\tau)$ from Eq.~(\ref{quadtau}), so that
\begin{equation} \label{sphi}
\phi=\phi_0+\frac{1}{n}\sum_{i\ge0}\varepsilon^i\phi^{(i)}(\tau).
\end{equation}
\par

Up to the second order of $\varepsilon$ we get
\begin{eqnarray} \label{qL}
q' &=& q^{(0)}(\tau)+q^{(1)}(\tau)+q^{(2)}(\tau), \\ \label{QL}
Q' &=& Q^{(0)}(\tau)+Q^{(1)}(\tau)+Q^{(2)}(\tau),
\end{eqnarray}
whereas Eq.~(\ref{sphi}) is better written in the fashion of Eq.~(\ref{fi4th}). Namely,
\begin{equation} \label{fitau}
\phi'=\phi'_0+\frac{1}{n}\left\{\omega\left(1+\frac{\Omega^2}{\omega^2}d\right)\tau+\frac{\Omega}{\omega}\big[p(\tau)-p(0)\big]\right\},
\end{equation}
in which, now, the modification of the linear trend given by $d$ and the long-period variation given by $p$, are computed in the Lindstedt series approach as
\begin{eqnarray} \label{dL}
d &=& d^{(0)}+d^{(1)}+d^{(2)}, \\ \label{pL}
p &=& p^{(0)}(\tau)+p^{(1)}(\tau)+p^{(2)}(\tau).
\end{eqnarray}
\par

Detailed expressions of the Lindstedt series solution were first presented in \cite{Lara2018ICATT}, and are reproduced in \ref{a:tables} for completeness.
\par

The orbital period of the long-term solution is then defined as
\begin{equation} \label{torbital}
T_O=\frac{2\pi}{\omega\left[1+(\Omega/\omega)^2d\right]},
\end{equation}
while the libration period is
\begin{equation} \label{tlibra}
T_L=\frac{2\pi}{\Omega{n}}.
\end{equation}
A periodic solution, on average, will happen when both periods are commensurable.
\par

We remark that the zeroth order terms as well as the terms $q^{(1)}$ and $Q^{(1)}$ of the first order of the Lindstedt series solution are equivalent to those computed in \cite{Lara2018} ---the term $\phi^{(1)}$ was not computed there, because the main concern was the correct modeling of the motion of the center of the reference ellipse, as opposite to the computation of precise ephemeris. Note, however, that, since we are using here a different scaling of the transformation to epicyclic variables, and we further made some rearrangement of the different orders of the Lindsdet series solution for computational efficiency, a direct comparison with the results in \cite{Lara2018} may not match with the lower orders of the solution provided here even if the new $k$ scaling is taken into account.

\section{Example applications}

The new Lindstedt series solution has been tested for different DRO. In particular, we found that it provides results of comparable accuracy to those obtained with the semi-analytical integration of the flow stemming from the long-term Hamiltonian (\ref{H8a}). Several cases are discussed below, which include the use of the orbit design parameters in the computation of periodic orbits.
\par

\subsection{Case 1. DRO with large libration}

\begin{figure}
\centering
\includegraphics[scale=0.88]{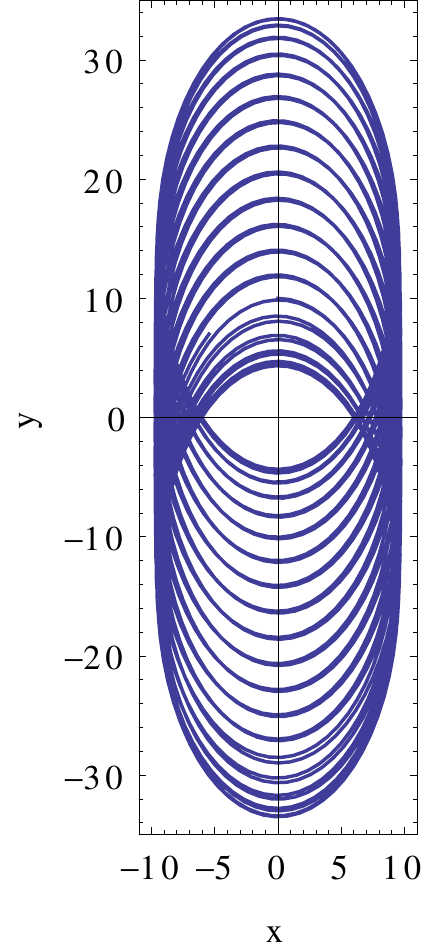}
\includegraphics[scale=0.88]{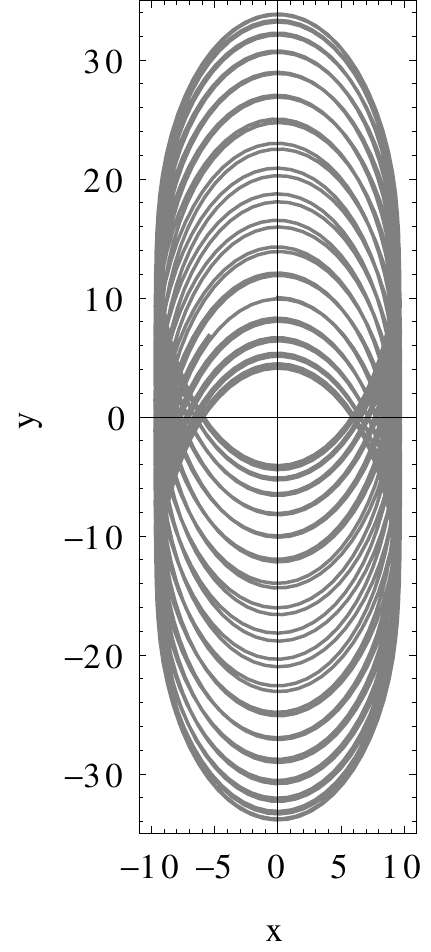}
\caption{Case 1. Left: numeric. Right: analytic.}
\label{f:orbit1}
\end{figure}

The first example presented is the case of orbits with large librational motion described in \cite{Lara2018}. Initial conditions for this sample orbit are $x=0$, $y=10$, $X=-0.5$, and $Y=-0.1$ in Hill units ($\omega=\mu=1$). For these initial conditions, the Lindstedt series solution predicts a libration frequency $\Omega=0.0187357$, less than 50 times smaller than the rotation rate of the system, an orbital period $T_O=6.27588$, and a libration period $T_L=231.669$. The numerical and analytical propagation of these initial conditions along one libration period are depicted in Fig.~\ref{f:orbit1}, which shows that the high-order Lindsted series solution now captures the bulk of the Hill problem dynamics of this kind of orbit.
\par

\begin{figure}
\centering
\includegraphics[scale=0.7]{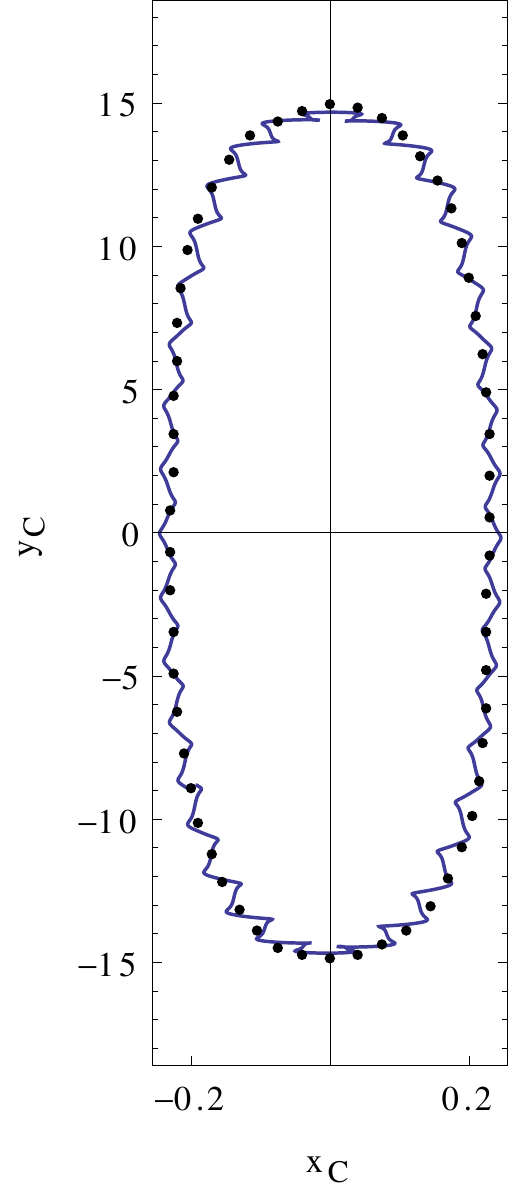}
\caption{Case 1. Evolution of the center of the reference ellipse. Dots: analytical solution; Full line: true (numerical) solution.}
\label{f:o1cel}
\end{figure}

Next, we check that the orbit described by the center of the reference ellipse is notably larger in the $y$ axis direction, exactly as predicted by the analytical solution. Indeed, as shown in Fig.~\ref{f:o1cel} the elongation of the center of the reference ellipse along the $x$ direction reaches only about $\pm0.2$ units, whereas the amplitude in the $y$ axis direction is about 75 times longer. Besides, the trajectory described by the prime variables averages quite accurately to the values provided by the true solution, as also shown in Fig.~\ref{f:o1cel}, where the true trajectory is displayed with a full line, and the prime variables analytical, long-term prediction is represented by dots.
\par

\begin{figure}
\centering
\includegraphics[scale=0.47]{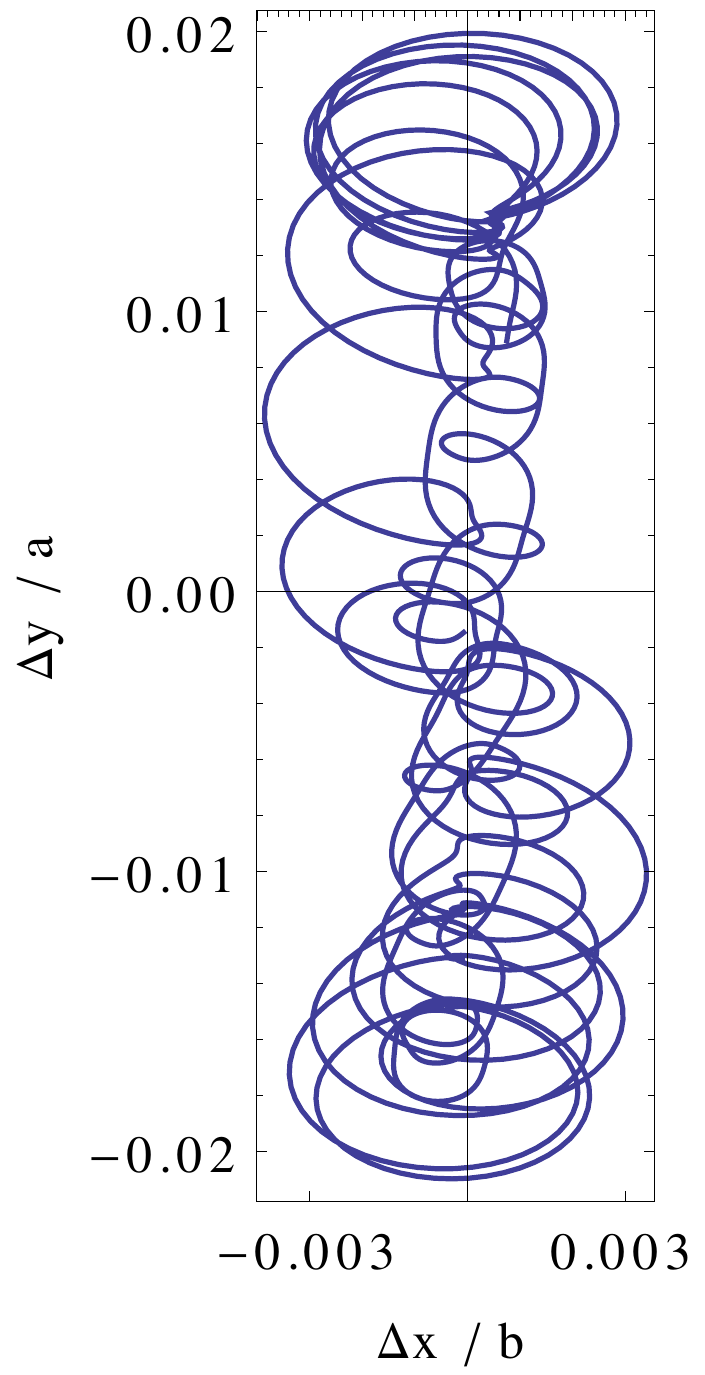}
\caption{Case 1. Scaled error in distance.}
\label{f:o1erdis}
\end{figure}

In addition to visually check the similitude of the analytical and numerical orbits in Fig.~\ref{f:orbit1}, the combined errors of the orbit in Cartesian coordinates are depicted in Fig.~\ref{f:o1erdis}, which presents the scaled error in distance of the true orbit with corresponding ephemerides computed from the  analytical solution in prime variables. To be more descriptive, the values of the errors in the $x$ and $y$ coordinates are scaled by the value of the semi-axis of the reference orbit in the respective direction, namely the errors in $x$ are scaled by $b$ and the errors in $y$ by $a=2b$. We calculated $\Phi'=45.1237$ for this orbit, which remains constan in the prime variables. Hence, we used the constant value $b=9.49987$, as computed from Eq.~(\ref{bxy}). This representation clearly shows that the errors in the $y$ axis direction are notable larger than those in the $x$ axis direction, which is expected from the longer excursion of the orbit along the $y$ axis.
\par

\begin{figure}[tb]
\centering
\includegraphics[scale=0.75]{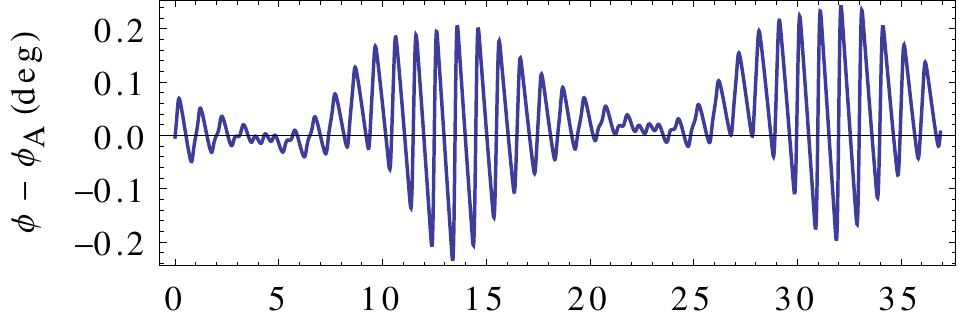} \\
\includegraphics[scale=0.75]{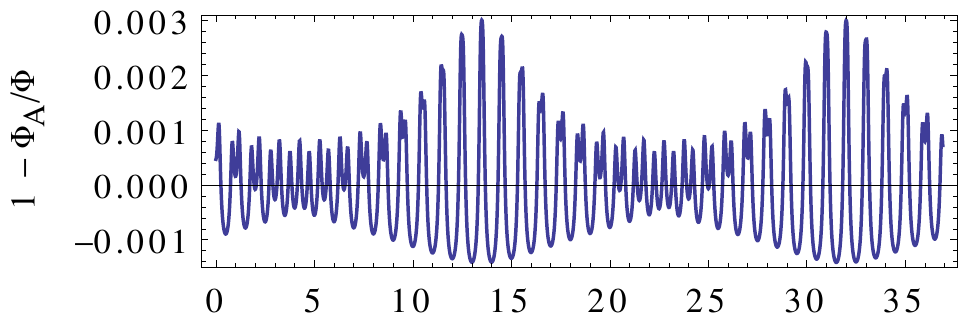} \\
\includegraphics[scale=0.75]{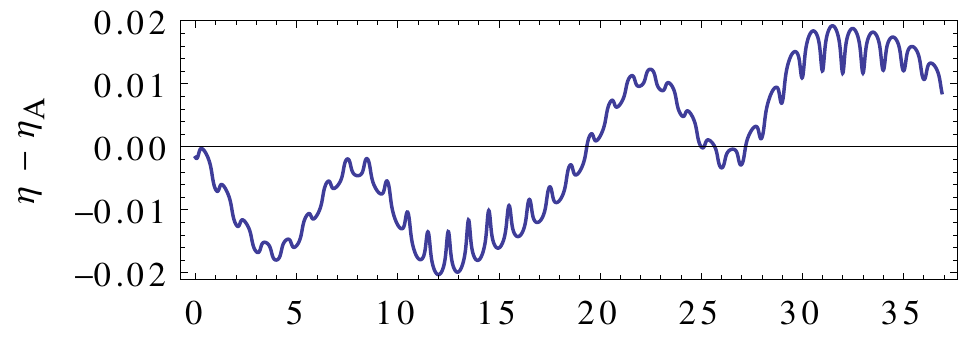} \\
\includegraphics[scale=0.75]{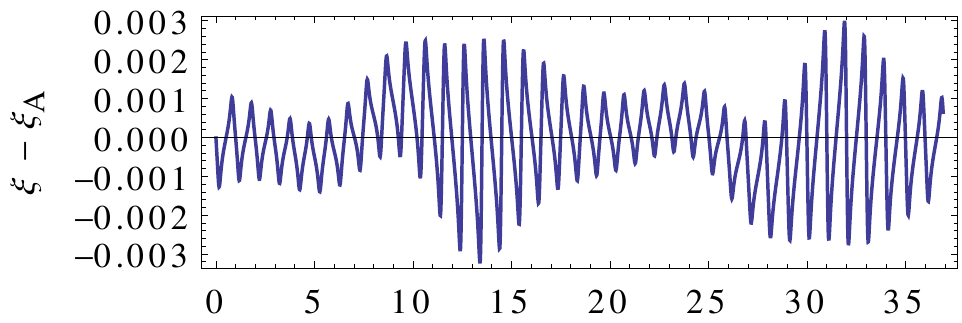}
\caption{Case 1. Errors in epicyclic variables between the true, numerical propagation and the Lindstedt series solution. Abscissas are orbital periods.}
\label{f:o1errepi}
\end{figure}

\begin{figure}[htb]
\centering
\includegraphics[scale=0.75]{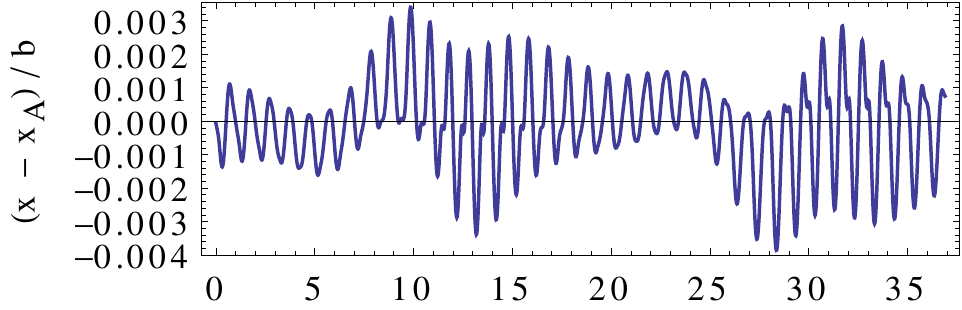} \\
\includegraphics[scale=0.75]{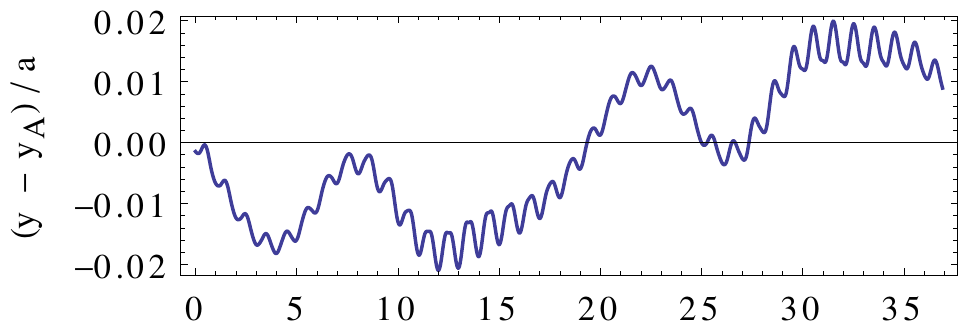} \\
\includegraphics[scale=0.75]{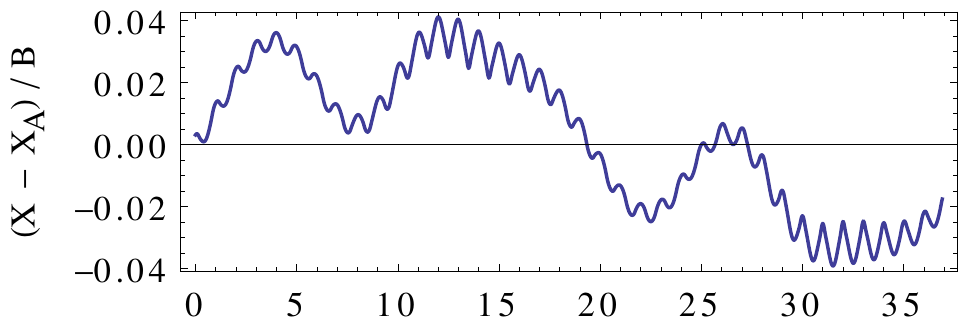} \\
\includegraphics[scale=0.75]{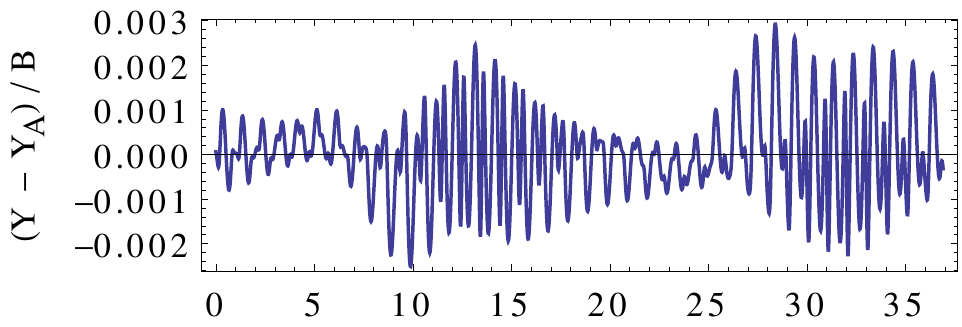}
\caption{Case 1. Scaled errors in Cartesian variables.}
\label{f:o1errxy}
\end{figure}

Higher detail on the errors between the analytical solution and the numerical propagation is provided in Fig.~\ref{f:o1errepi}, where the time history of the errors is displayed in epicyclic variables, and in Fig.~\ref{f:o1errxy}, that depicts the errors in Cartesian variables. The latter are obtained by the direct transformation of the prime epicyclic variables to prime Cartesian using Eqs.~(\ref{tdir})--(\ref{bxy}).
\par

The errors of the analytical propagation comprise both, short- and long- period effects, of notable amplitude, as well as a small secular trend. Thus, we observe in the top plot of Fig.~\ref{f:o1errepi} that the error in the phase of the orbit oscillates along one orbital period with an amplitude of just a few tenths of degree with respect to the true value. These short-period oscillations seem to be modulated with long-period oscillations of half the libration period. A more detailed examination of this plot reveals a small secular trend of just a few hundredths of degree by libration period, which remains almost hidden by the large amplitude of the periodic oscillations.
\par


The relative errors between the true value of $\Phi$ and the constant one $\Phi'$ remain of the order of one thousandth, as shown in 2nd plot of Fig.~\ref{f:o1errepi}. These errors are mostly a consequence of comparing prime with original variables, without having recovered the short-period effects that were removed by the perturbation theory.
\par

An analogous propagation of the errors of the other epicyclic elements is observed in the third and bottom plots of Fig.~\ref{f:o1errepi}, showing again the short and long period modulation. The worst case happens to $\eta$ ($q$ variable), where the amplitude of the long-period oscillation of the errors is notably larger than the amplitude of the short-period oscillation. This behavior is as a consequence of the early truncation of the perturbation theory for an orbit with high amplitude oscillations in this variable.
\par

As expected from the transformation equations to epicyclic variables given in Eqs.~(\ref{tdir})--(\ref{bxy}), the errors in $q$ propagate directly to the $y$ coordinate and the $X$ conjugate momentum, as can be observed in the second and third plots of Fig.~\ref{f:o1errxy}. Note that the errors of $X$ and $Y$ are also scaled, in this case by the radius of the reference circumference of the hodograph $B$.
\par

The errors of $x$ remain of the order of just a few thousandths when compared with the orbit dimension in the $x$ axis direction, and the same happens to the error in $Y$ when compared to the size of the hodograph, as shown in the top and bottom plots of Fig.~\ref{f:o1errxy}, respectively. On the contrary, the errors of $y$ are one order of magnitude higher when compared with the orbit dimension in the $y$ axis direction, and the same increase in the magnitude of the errors happens to $X$ relative to the size of the hodograph, as illustrated with the second and third plots of Fig.~\ref{f:o1errxy}, respectively.
\par

\subsection{Case 2. Short-period effects}

\begin{figure}[tb]
\centering
\includegraphics[scale=0.7]{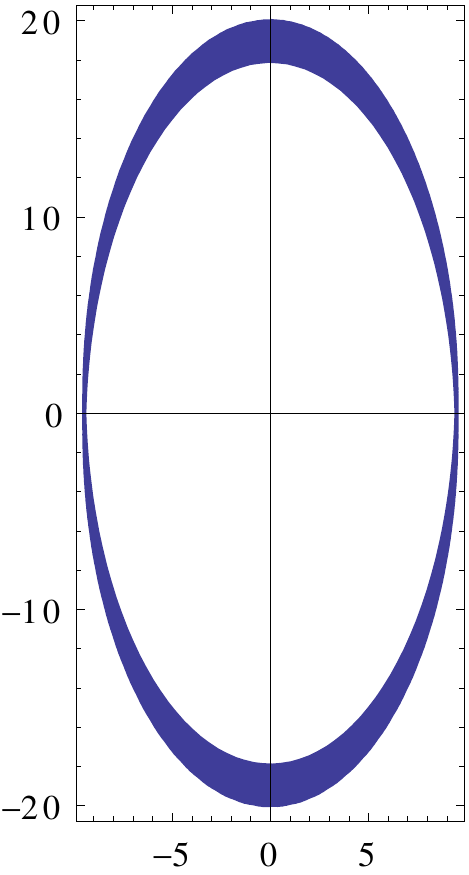} \quad
\includegraphics[scale=0.7]{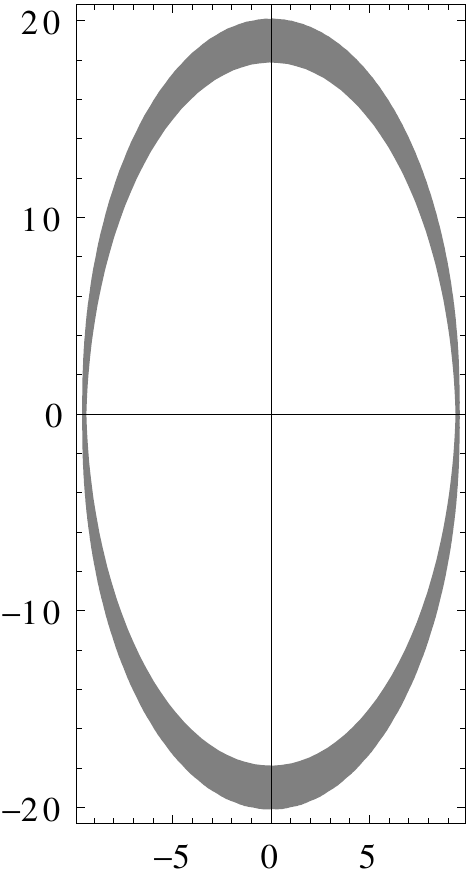}
\caption{Case 2. Left: numeric. Right: analytic.}
\label{f:orbit2}
\end{figure}

The second example is the the test case 2 in \cite{Lara2018} for small amplitude librations, which has been chosen to illustrate the effect of the short-period corrections. The initial conditions are now $x=0.1$, $y=20.0$, $X=-10.5$, and $Y=-0.1$, for which the higher order Lindstedt series solution predicts an orbital period $T_O=6.27815$ and a libration period $T_L=334.835$. The true, numerical solution is propagated for one libration period and displayed in the left plot of Fig.~\ref{f:orbit2}, where the analytical prediction is also displayed for the same period in the right plot. As shown in the figure, both orbits match at the precision of the graphics.
\par

\begin{figure}[h]
\centering
\includegraphics[scale=0.7]{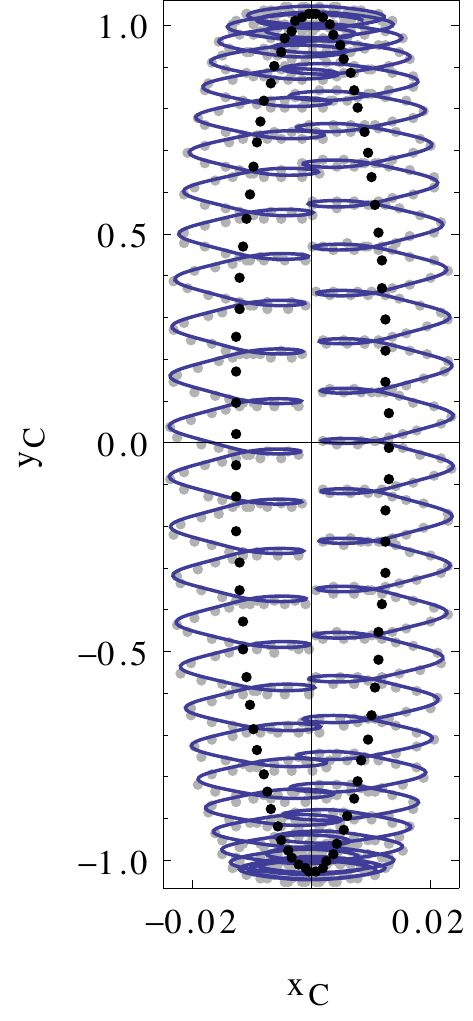}
\caption{Case 2. Analytic (dots) and numeric (full line) evolution of the center of the reference ellipse. The gray dots include short-period corrections.}
\label{f:o2cel}
\end{figure}

The center of the reference of the ellipse describes now a more intricate path, whose evolution is depicted in Fig.~\ref{f:o2cel} for one librational period. The full line corresponds to the numerical solution, while dots correspond to analytical predictions. The black dots show that, on average, the center of the reference ellipse moves also on an ellipse, as predicted by the theory in prime variables. Besides, we see that the results predicted by the analytical solution, represented with gray dots in Fig.~\ref{f:o2cel}, superimpose to the true trajectory when the short-period corrections are applied to recover the original variables.
\par

\begin{figure}[t]
\centering
\includegraphics[scale=0.7]{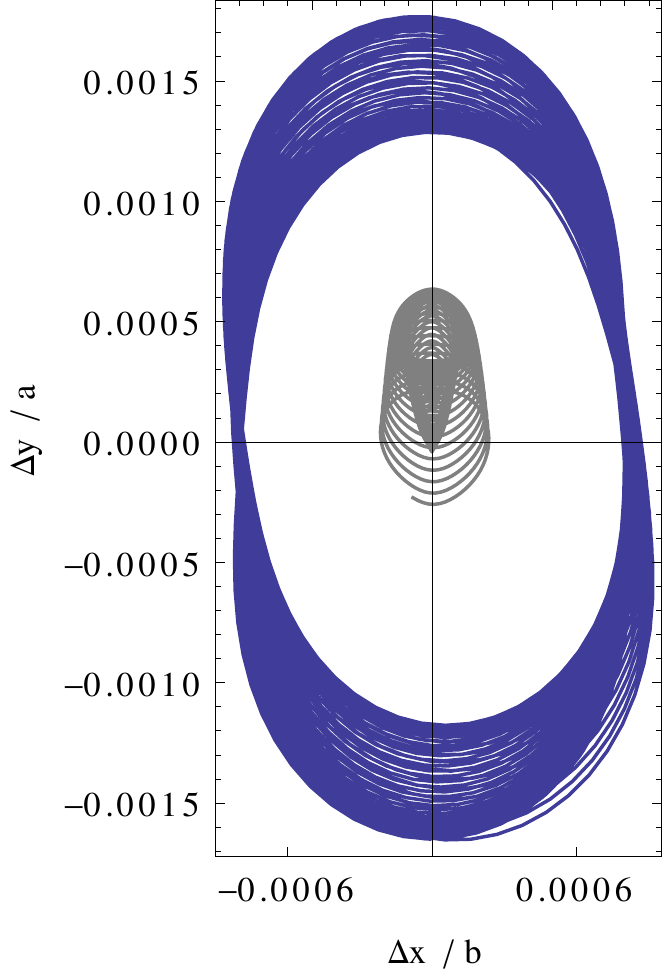}
\caption{Case 2. Scaled error in distance for a libration period when the true orbit is compared with prime, mean elements or osculating ones (gray curve) obtained from the analytical solution.}
\label{f:o2erdis}
\end{figure}

The improvements obtained when recovering the short-period effects of the analytical solution are better appreciated in Fig.~\ref{f:o2erdis}, where the magnitude of the errors of the analytical orbit, in the previously mentioned scaling by the length of each semi axis of the reference ellipse, are displayed in blue for the prime variables, and in gray for the full solution. The fact that the prime variables propagation starts from non-null error, is due to the fact that we are comparing sets of variables of different nature. This also happens to the errors displayed in gray, yet to a very minor extent. In this last case, the variables compared enjoy the same nature, and the initial error is a consequence of the truncation order of the perturbation theory, which makes that going forth with the inverse transformation, and back with the direct one, does not give exactly the same initial conditions.
\par

\begin{figure}[htb]
\centering
\includegraphics[scale=0.75]{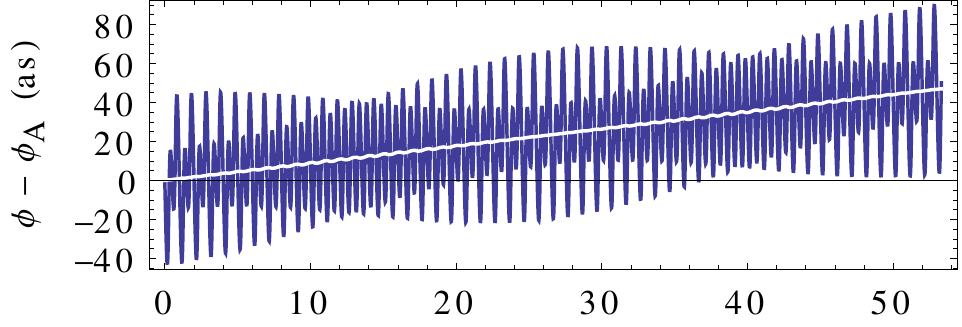} \\
\includegraphics[scale=0.75]{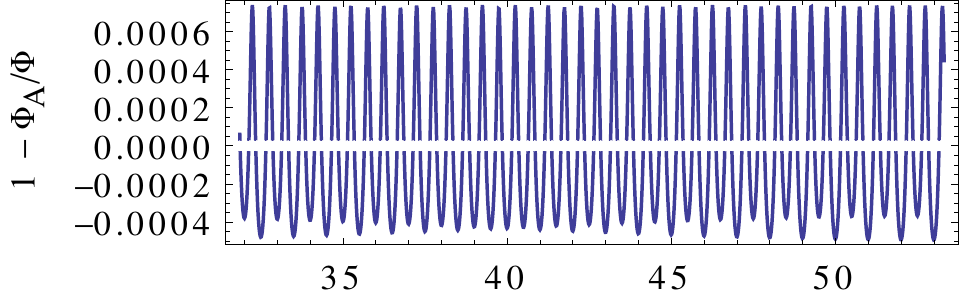} \\
\includegraphics[scale=0.75]{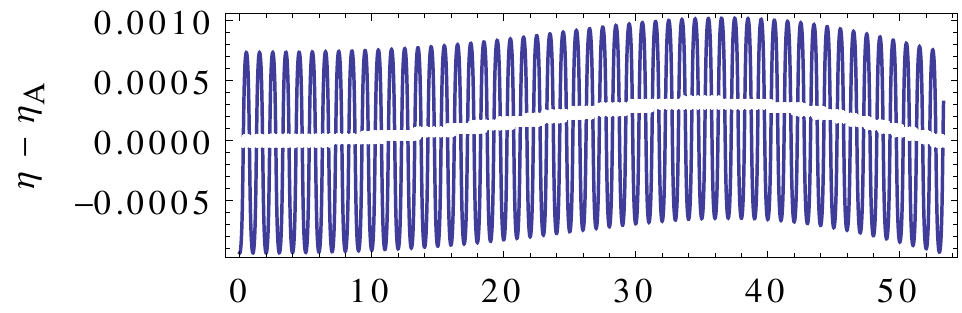} \\
\includegraphics[scale=0.75]{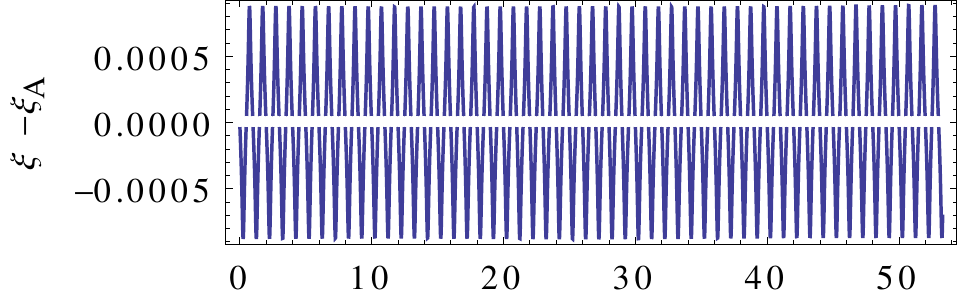} 
\caption{Case 2 orbit. Errors of the mean elements propagation in epicyclic variables with the errors after recovering the short-period effects superimposed in white. Abscissas are orbital periods.}
\label{f:o2errepi}
\end{figure}

Figure \ref{f:o2errepi} depicts the time history of the errors in epicyclic variables. Short-period effects (in blue, in the figure) are mostly flattened when the direct transformation is applied to the prime epicyclic variables (gray curves). Still, a long-period modulation is observed in the time history of the errors of $\eta$ (variable $q$). Besides the errors of the phase of the satellite predicted by the analytical solution are affected of a secular increase at a linear rate of about $0.88$ as times orbital period.
\par

\begin{figure}
\centering
\includegraphics[scale=0.75]{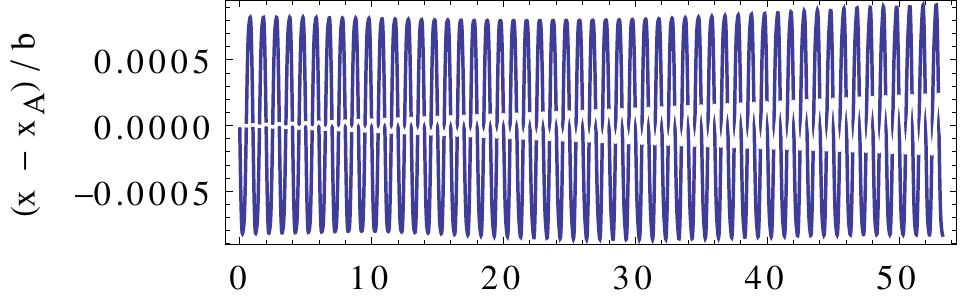} \\
\includegraphics[scale=0.75]{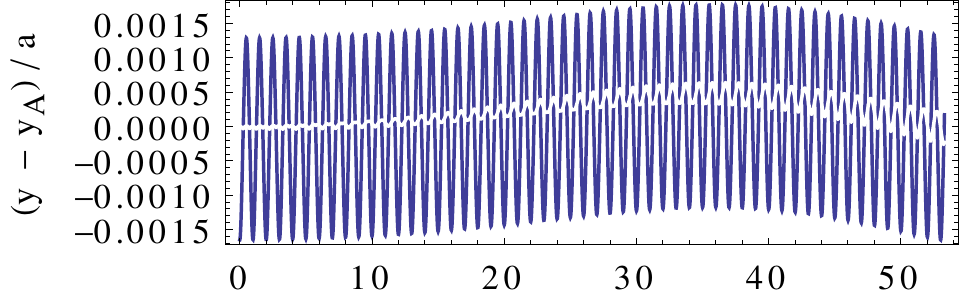}
\includegraphics[scale=0.75]{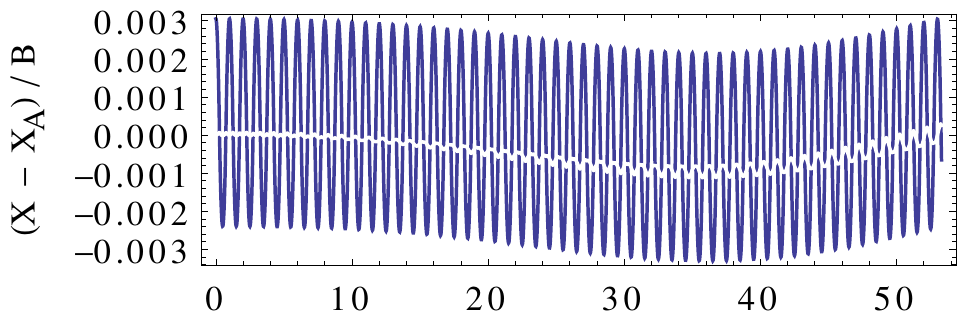}
\includegraphics[scale=0.75]{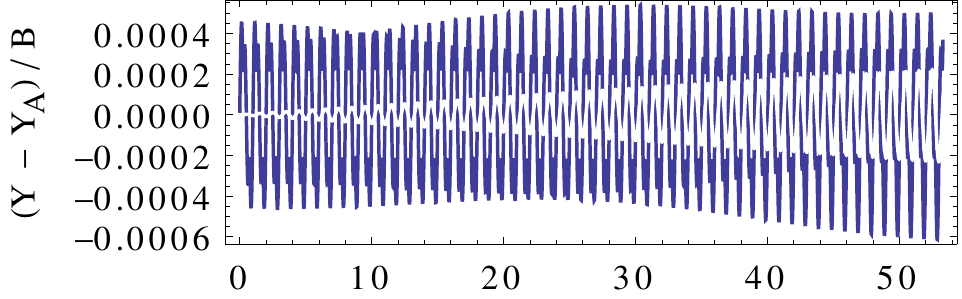}
\caption{Case 2 orbit. Errors of the mean elements propagation in Cartesian variables with the errors after recovering the short-period effects superimposed in white.}
\label{f:o2errxy}
\end{figure}

The propagation of the errors in Cartesian variables is shown in Fig.~\ref{f:o2errxy}; again, these errors are scaled to the size of the reference orbit. Now, while applying the periodic corrections clearly improve the precision of the analytical solution, the secular trend observed in the propagation of $\phi'$, which is due to the early truncation of the perturbation theory, affects to the error propagation of all the Cartesian variables. In any case, the errors remain in the order of one thousandth relative to the orbit, even when the short-period effects are not corrected.
\par

\subsection{Case 3. Close approach}

\begin{figure}
\centering
\includegraphics[scale=0.7]{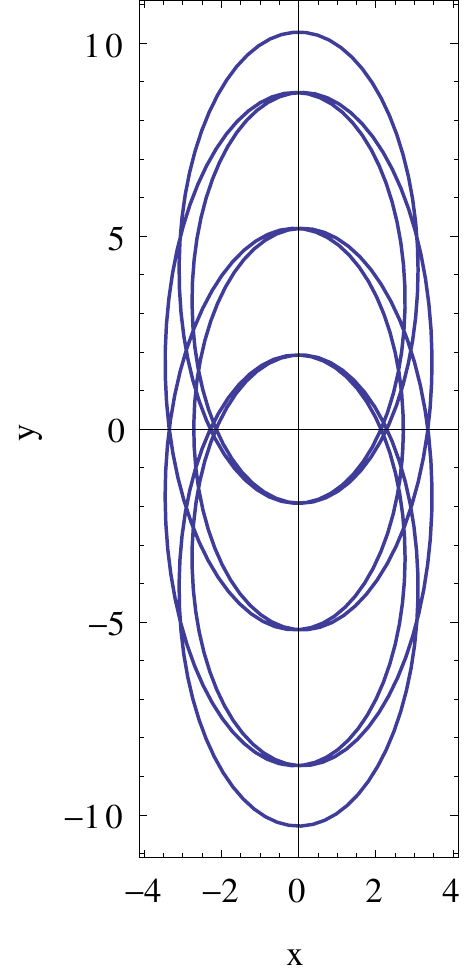}\qquad
\includegraphics[scale=0.7]{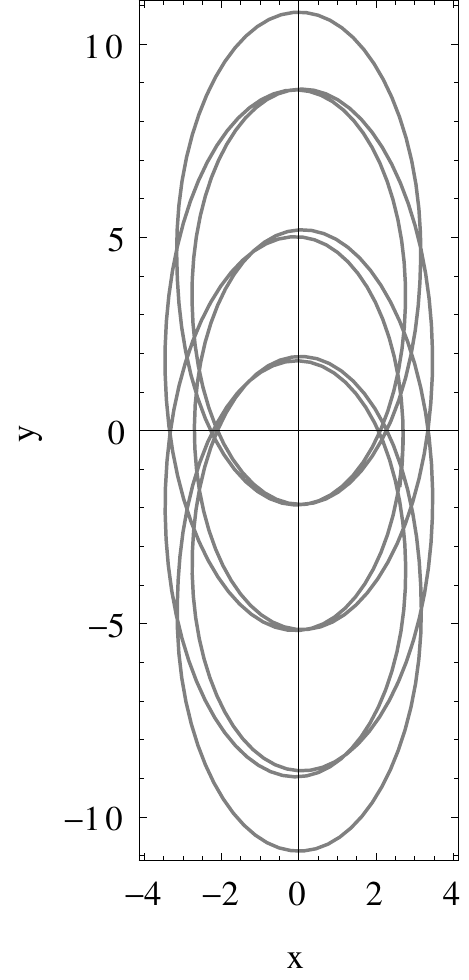}
\caption{Case 3. Left: numerical propagation. Right: analytical solution.}
\label{f:orbit3}
\end{figure}

To test the limits for the validity of the new analytical solution, we choose a extreme case which is borrowed from \cite{LidovVashkovyak1994a}. Initial conditions for this orbit are $x=2.7163$, $y=0$, $X=0$, $Y=-2.9724$, and correspond to an orbit that can approach to the primary as close as $\sim2$ Hill units in the $y$ axis direction, and evolves quite fast, with a libration period that is only slightly longer than 7 times the orbital period.  The numerical and analytical propagation of these initial conditions for one librational period are depicted in Fig.~\ref{f:orbit3}. While quantitative variations are clearly apparent in that figure, it can be observed that the analytical solution predicts correctly the main characteristics of the orbit evolution.
\par

\begin{figure}
\centering
\includegraphics[scale=0.7]{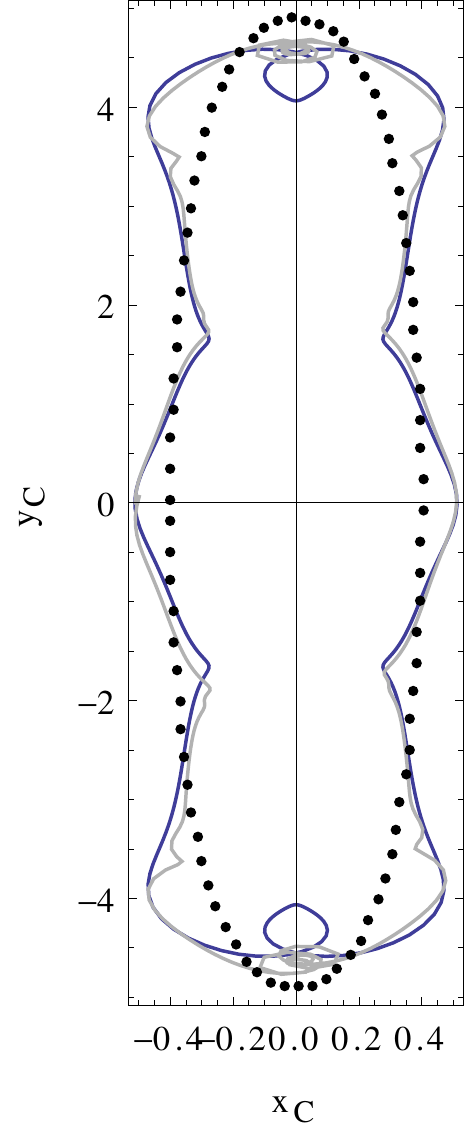}
\caption{Case 3 orbit. True, numeric evolution of the center of the reference ellipse with analytical mean (dots) and osculating predictions (gray line) superimposed. }
\label{f:o3cel}
\end{figure}

In fact, in spite of the irregular behavior of the center of the reference ellipse of this orbit, with two loops at the extreme positions along the $y$ axis, the prime variables solution seems to correctly average the true trajectory, as depicted with the black dotted curve in Fig.~\ref{f:o3cel}. On the other hand, the full analytical solution (in gray line in Fig.~\ref{f:o3cel}), which includes the short-period effects of the perturbation solution, adheres quite closely to the real trajectory (the dark solid line in Fig.~\ref{f:o3cel}), even to the extent of detecting the loops at the $y$ axis, yet failing to reproducing them accurately.

\subsection{Case 4: single periodic DRO}

For the typical 1:1 DRO, we choose $a=\rho=10$ in Hill units, and compute $\Phi'=12.5$. Despite $\Omega=0.0490672$ for this value of $a$, and hence our solution predicts a libration period about 20 times higher than the orbital one, this would not fit to the 1:1 resonance. Indeed, if we fix $q'_0=0$, the fact that $a$ and $\rho$ are equal makes that $Q'_0=0$, which immediately shows that terms in the square brackets of Eqs.~(\ref{qpL}) and (\ref{QpL}) always vanish, resulting in a fixed guiding center at the origin, on average. Therefore, we impose $T_L=T_O=6.24852$, as obtained from Eq.~(\ref{torbital}).
\par

\begin{figure}[ht]
\begin{center}
\includegraphics[scale=0.75]{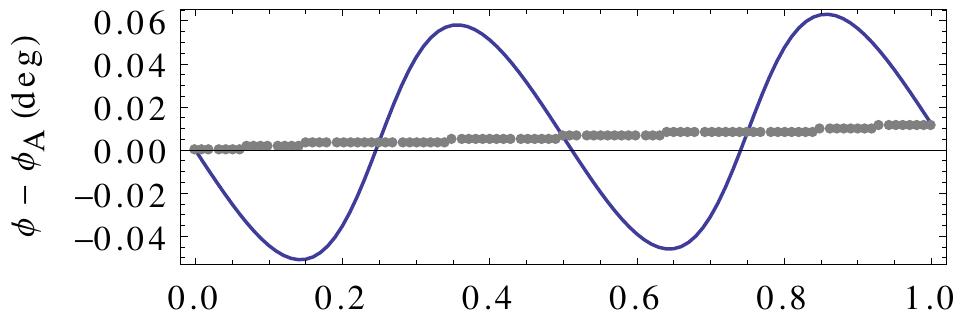} \\[1ex]
\includegraphics[scale=0.75]{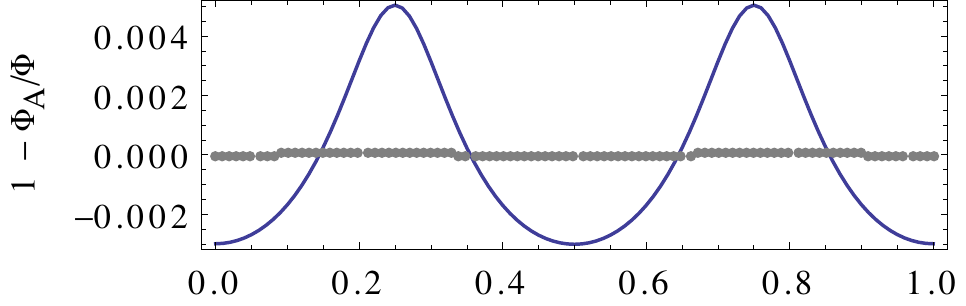} \\[1ex] 
\includegraphics[scale=0.75]{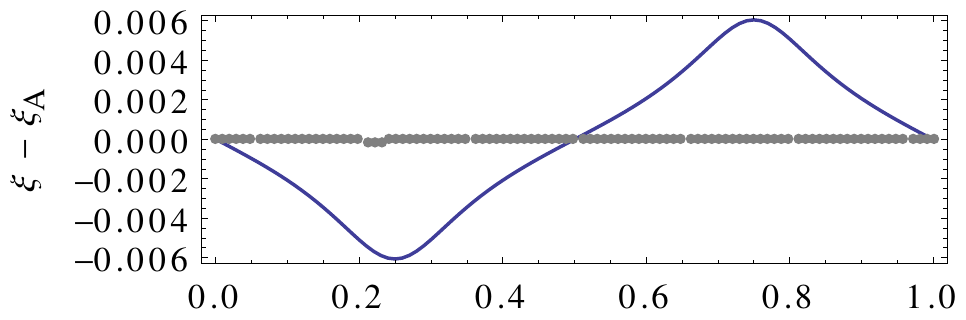} \\[1ex]
\includegraphics[scale=0.75]{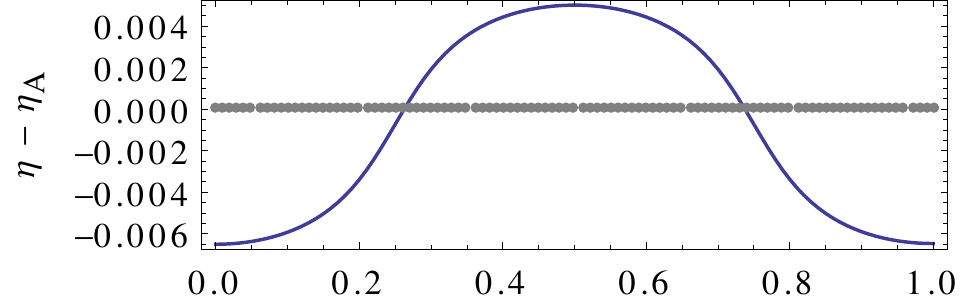}
\caption{Errors of the analytical orbit with design parameters $a=\rho=10$. Solid line: mean elements prediction. Gray dots: osculating elements prediction. Abscissas are orbital periods.}
\label{f:a10p10err}
\end{center}
\end{figure}

The time history of the errors of this analytical solution is shown in Fig.~\ref{f:a10p10err}. The errors of the mean elements solution are very small, and almost negligible when the short-period corrections are taken into account (light dots in Fig.~\ref{f:a10p10err}). The larger errors  occur for the phase of the orbiter $\phi$, where a secular trend due to the early truncation of the theory is clearly observed.
\par

If desired, the initial conditions provided by the analytical solution can be improved with differential corrections to get an exact periodic orbit in the 1:1 resonance. Using the algorithm of \cite{LaraPelaez2002}, we find that just 3 iterations are enough to obtain the initial conditions of a true periodic, stable, DRO of the Hill problem with initial conditions $x=Y=0$, $y=9.783444749944893$, and $X=-4.847560254601411$, that, when propagated for an improved period $T_O=6.247084797518564$ yields periodicity errors of the order of $10^{-12}$ for $x$ and $Y$, and $10^{-10}$ for $y$ and $X$.
\par

\subsection{Case 5. Multiple periodic DRO}

Another example is presented for a periodic DRO with larger libration, which is definitely much more challenging case. Again we choose $a=10$, but now with $\rho=5$. Next, we compute $\Phi'=12.5$ and $\Omega=0.0490672$, about 20 times smaller than $\omega$. Then, for $q'_0=0$ we compute $Q'_0=0.141645$ and choose $\phi'_0=0$. However, if we compute the orbital and librational periods from Eqs.~(\ref{torbital}) and (\ref{tlibra}), respectively, we check that $r=T_L/T_O=18.29$ showing that the orbit with design parameters $a=10$, $\rho=5$ is not periodic.
\par

The commensurability between the orbital and librational periods can be improved to get $r^*=18$ by a sequence of iterations that leaves fixed $\rho$ and $q'_0=0$, and modifies $\Phi'$, and, therefore, $\Omega$, $a$, and $Q'_0$. Indeed, in view of $r={T_L(q'_0,Q'_0,\Phi')}/{T_O(q'_0,Q'_0,\Phi')}$, the secant method can be used to compute $r_i=r_i(\Phi_i)$ such that the difference $r^*-r_i$ is as small as desired. Thus, after 3 iterations we find that choosing $a=9.87661$, $\rho=5$, yields the required periodicity after 18 orbits, which are traveled in a librational period $T_L=112.379$.
\par

\begin{figure}[t]
\begin{center}
\includegraphics[scale=0.75]{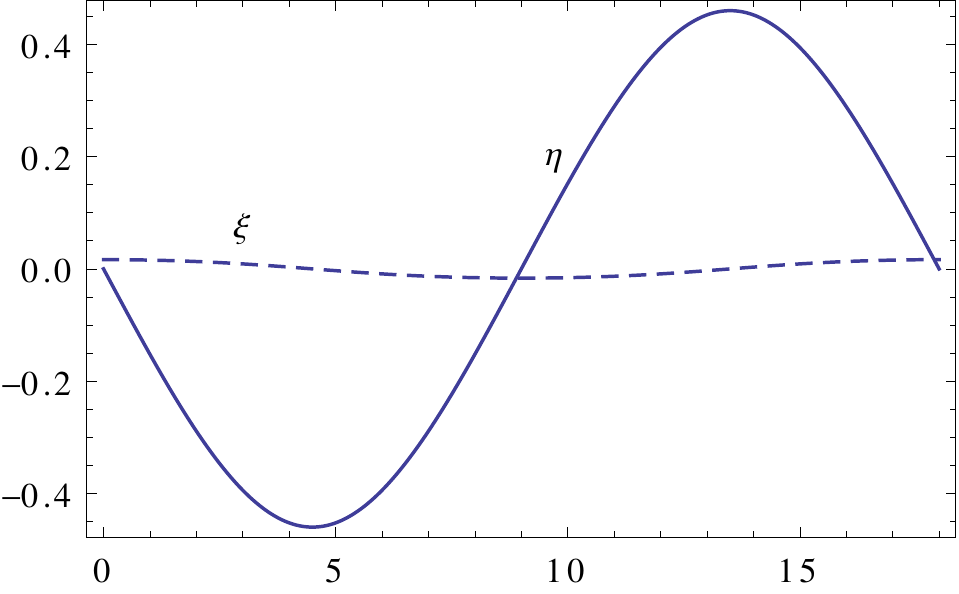} 
\caption{Non-dimensional epicyclic variables $\xi$ (dashed line) and $\eta$ (solid line) along one libration period of the orbit with design parameters $a=10$, $\rho=5$. Abscissas are orbital periods.}
\label{f:a10p5xieta}
\end{center}
\end{figure}

The analytical propagation of these initial conditions for a libration period shows that $\phi'$ grows with an almost linear rate which is just slightly higher than $2\pi$ times orbital period; $\Phi'$, of course, remains constant, and the evolution of the non-dimensional version of $Q'$ and $q'$ given by $\eta$ and $\xi$, respectively, is depicted in Fig.~\ref{f:a10p5xieta}.

\begin{figure}[h]
\begin{center}
\includegraphics[scale=0.7]{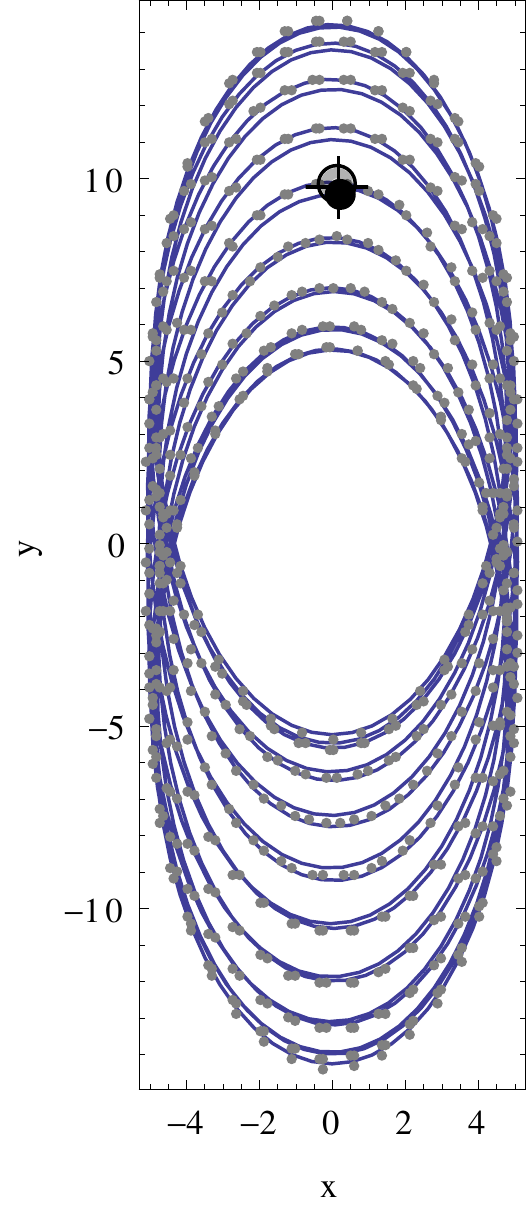} \quad
\includegraphics[scale=0.7]{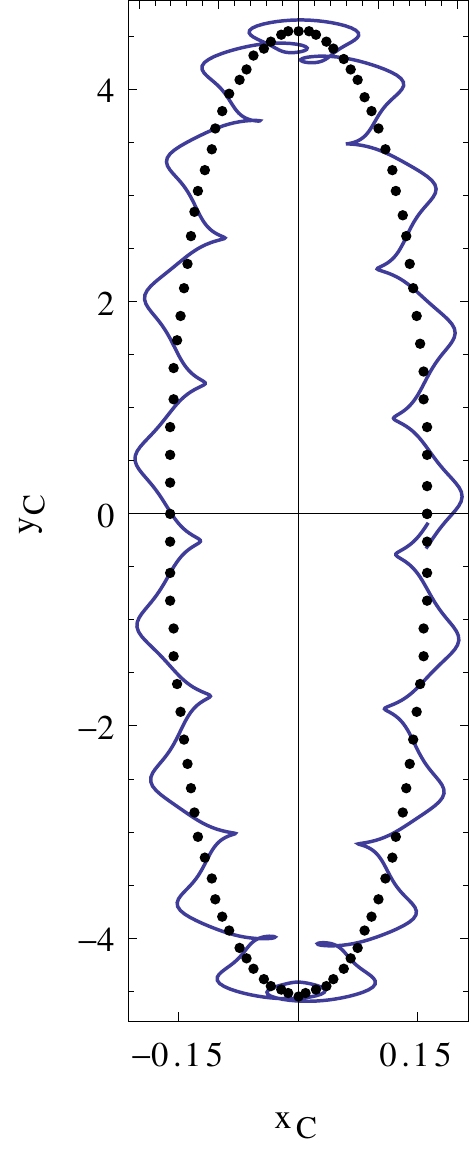} 
\caption{Orbit with design parameters $a=9.87661$, $\rho=5$, evaluated along one libration period $\sim18$ orbital periods. Left: Mean orbit prediction (gray dots) superimposed to the numerical propagation (solid line); the big black dot at $\sim(0.19,9.5)$ marks the last point for the predicted libration period. Right: Mean evolution of the guiding center (black dots) superimposed to the true trajectory (solid line).}
\label{f:a10p5}
\end{center}
\end{figure}

If we now compute initial conditions in the original, non averaged space and compute the corresponding orbit by numerical propagating the equations of motion of the Hill problem, we find that the real orbit and the analytical prediction match quite well, as depicted in the left plot of Fig.~\ref{f:a10p5} where the mean orbit is represented with light dots and the real orbit with a solid line. On the other hand, the predicted behavior of the guiding center indeed averages the trajectory derived from the real trajectory, as shown in the right plot of Fig.~\ref{f:a10p5}, where the analytical prediction is represented with dots.
\par

\begin{figure}[htb]
\begin{center}
\includegraphics[scale=0.75]{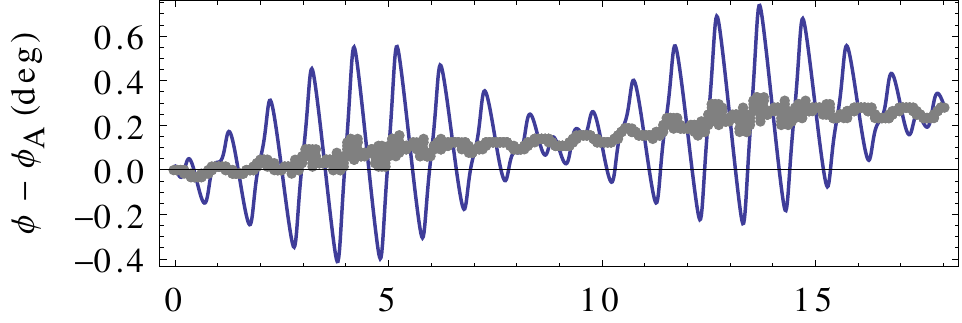} \\[1ex]
\includegraphics[scale=0.75]{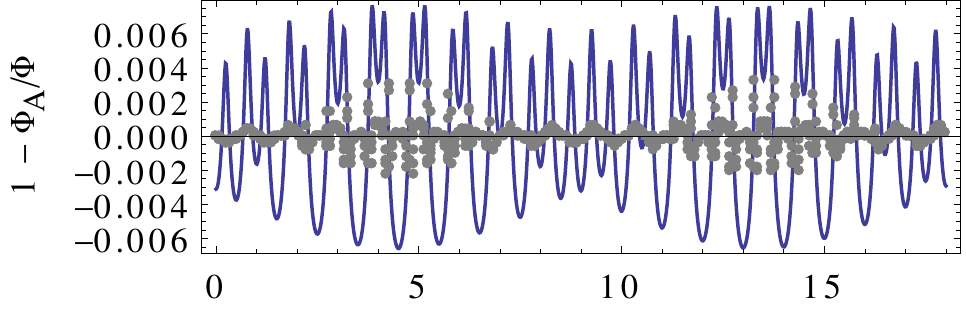} \\[1ex] 
\includegraphics[scale=0.75]{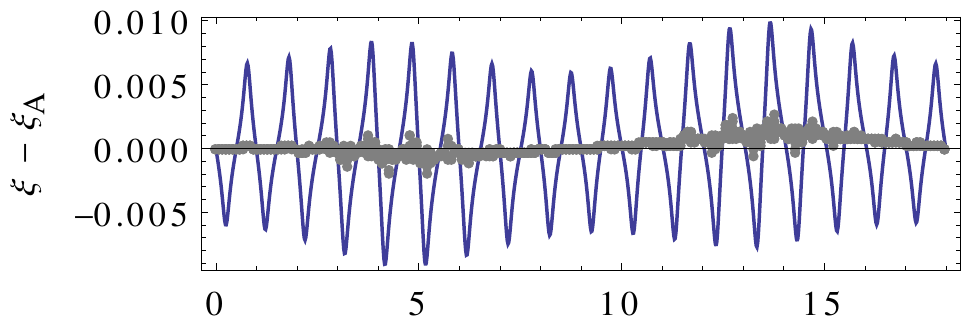} \\[1ex]
\includegraphics[scale=0.75]{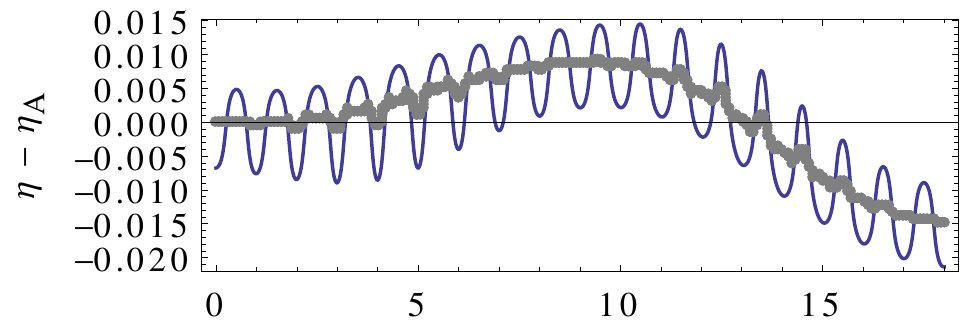}
\caption{Errors of the orbit with design parameters $a=9.87661$, $\rho=5$, with respect to the mean elements (solid line), and osculating elements prediction (gray dots). Abscissas are orbital periods.}
\label{f:a10p5err}
\end{center}
\end{figure}

Things fit even much better when short periodic effects are added to the mean elements prediction. Indeed, as shown in Fig.~\ref{f:a10p5err}, short-period effects are the more important source of errors of the mean analytical solution (solid lines), whereas osculating predictions flatten the errors (light dots) and clearly disclose the secular and long-period errors that result from the early truncation of the perturbation approach.
\par

\begin{figure}[htb]
\begin{center}
\includegraphics[scale=0.7]{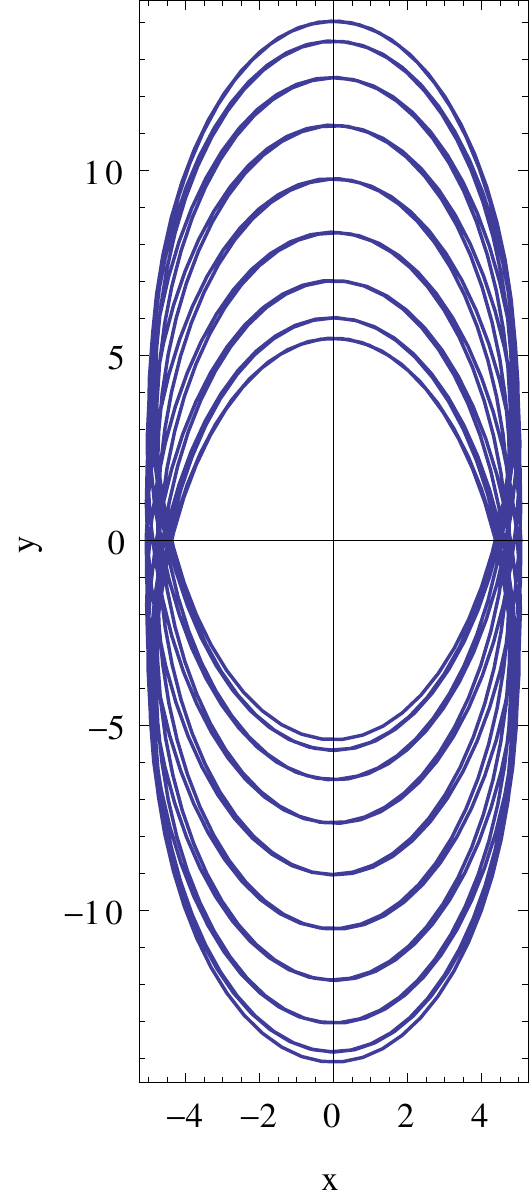} \quad
\includegraphics[scale=0.7]{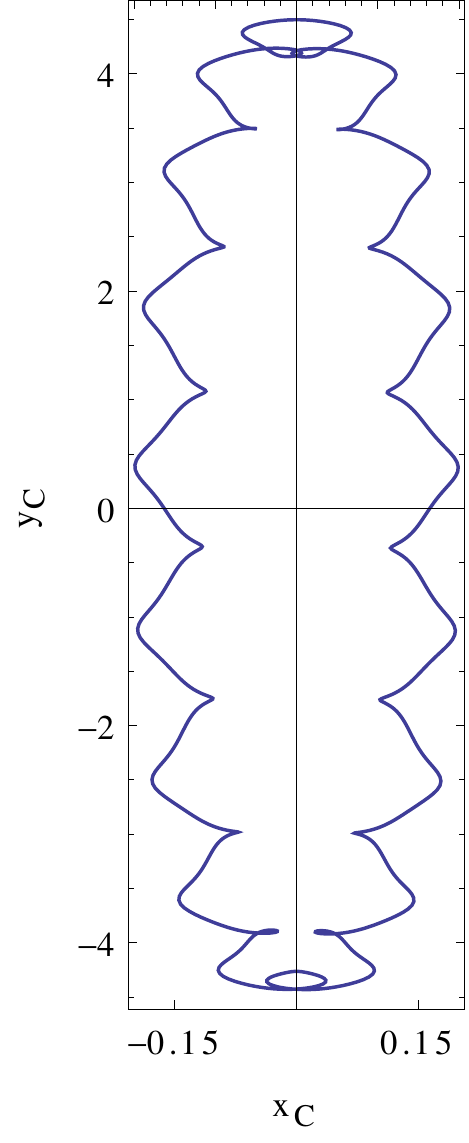}
\caption{True periodic orbit corresponding to the improved design parameters $a=9.87661$, $\rho=5$. Left: numerically integrated orbit. Right: center of the reference ellipse.}
\label{f:improved}
\end{center}
\end{figure}

If desired, the initial conditions and period provided by the analytical solution are further amenable of improvement by differential corrections to give a true periodic orbit. Application of differential corrections to the improved initial conditions and period with the algorithm of \cite{LaraPelaez2002}, needs just 4 iterations to converge to a true periodic, stable, DRO of the Hill problem with initial conditions 
\begin{eqnarray*}
x &=& 5.061558354876498, \\
y &=& 0, \\
X &=& 0.1831185556870679, \\
Y &=&-5.003556180647312,
\end{eqnarray*}
and period $112.3791870019849$, producing periodicity errors $\mathcal{O}(10^{-14})$ for ${x}$ and ${Y}$, and $\mathcal{O}(10^{-13})$ for ${y}$ and ${X}$. The resulting true periodic orbit is illustrated in Fig.~\ref{f:improved}.
\par

\section{Conclusions}

Useful analytical solutions for the preliminary design of distant retrograde orbits are efficiently computed based on Lidnstedt series. These kinds of perturbation solutions effectively capture the main features of the dynamics even in extreme cases in which the orbiter gets close to the co-orbital primary in its apparent journey about it. The solutions are time explicit and provide similar accuracy to the one achieved with previous approaches in the literature based on semi-analytical integration. The solution shows the existence of two design parameters, which are related to the size of the orbit and the minimum libration distance to the co-orbital primary. The proper choice of these parameter provides the initial conditions of an orbit with the characteristics determined by the user that is periodic on average. These initial conditions are amenable to the standard improvement with differential corrections, so that a true periodic orbit of the Hill problem with the desired attributes is easily obtained.
\color{black}

The Hill problem has been chosen as a demonstration problem because of its simplicity and generality. However, the approach taken here can be easily extended to more realistic dynamical models, like the restricted three-body problem or any of its variations resulting when different perturbations are taken  into account.

\subsection*{Declaration of interest}

None

\subsection*{Acknowledgemnts}

Support by the Spanish State Research Agency and the European Regional Development Fund under Projects ESP2016-76585-R and ESP2017-87271-P (MINECO/ AEI/ERDF, EU) is recognized.

\appendix

\section{Periodic corrections} \label{a:peco}

For completeness, we provide the explicit expressions of the transformation of the short-period elimination up to the truncation order, which are essentially the same as those in \cite{Lara2018} except for the different $k$ scale. Besides, the $\arctan(2\cot\phi)$ that appeared in the 8th order correction to the phase $\phi$ ---which without the proper modulation constrains this correction to the interval $(-\frac{\pi}{2},\frac{\pi}{2})$--- is now replaced by half the elliptic integral of the third kind of parameter $3/4$, argument $\phi$, and null elliptic modulus, thus allowing the phase to grow unbounded.

We recall that, for conciseness, we abbreviated $c\equiv\cos\phi$, $s\equiv\sin\phi$, and defined the auxiliary quantities $\tilde{K}=K(k^2)/\pi$, $\tilde{E}=E(k^2)/\pi$. Besides, we define the functions
\begin{eqnarray}
F^* &=& 2\tilde{K}\,\phi-F(\phi\,|\,k^2) \\
E^* &=& 2\tilde{E}\,\phi-E(\phi\,|\,k^2) \\
P^* &=& 2\phi-\Pi(k^2;\phi\,|\,0)
\end{eqnarray}
where $F(\phi\,|\,k^2)$ and $E(\phi\,|\,k^2)$ note the incomplete elliptic integral of first and second kind, respectively, and $\Pi(\alpha^2;\phi\,|\,\kappa^2)$ is the elliptic integral of the third kind. Note that $F^*$, $E^*$, and $P^*$, are periodic functions of $\phi$ with period $\pi$.

Then, the motion in the original coordinates, up to the truncation order of the perturbation approach, is obtained by adding the short-period effects to the epicyclic prime coordinates. Thus, for a generic element $\beta\in(\phi,q,\Phi,Q)$, the direct transformation is
\[
\beta=\beta'+\gamma(\Phi')\sum_{1\le{i}\le6}\delta\beta_i(\phi',q',\Phi',Q'),
\]
where the non-vanishing corrections are
\begin{eqnarray*}
\delta\phi_1 &=& -\frac{1}{2}F^* 
\\
\delta\phi_2 &=& -\frac{\eta}{\Delta}s
\\
\delta\phi_3 &=& 
\frac{c}{\Delta}\Big[\Big(\frac{1}{\Delta ^2}+1\Big) \eta ^2 k^2 s+2\xi \Big] +\eta^2[E^*-F^*] 
\\
\delta\phi_4 &=& \eta\Big[\Big(\frac{1}{\Delta ^3}-8 \tilde{E}\Big) \xi
+\frac{1}{9}\frac{s}{\Delta}\Big(\frac{3}{\Delta ^4}-\frac{7}{\Delta ^2}-14\Big)\eta^2\Big] 
\\
\delta\phi_5 &=& \gamma\Big[P^*-\Big(2k^2\tilde{K}+\frac{1}{4\Delta}\Big)F^*\Big] \\
&&+\frac{5}{36}(14E^*-11F^*)\eta^4+(F^*-4E^*)\xi^2 \\
&& +\frac{c}{\Delta}\Big[\frac{2}{3}\Big(8+\frac{1}{\Delta^2}-\frac{3}{\Delta^4}\Big)\eta^2\xi 
-\Big(3+\frac{k^2}{\Delta^2}\Big)\xi^2s \\
&& +\frac{5}{48}\Big(14+\frac{11}{\Delta^2}+\frac{8}{\Delta^4}-\frac{5}{\Delta^6}\Big)\eta^4s \Big],
\\
\delta\phi_6 &=& \frac{\gamma\eta}{\Delta}
\Big[ s\Big(\frac{11}{8 \Delta }-4\tilde{K}\Big)+\frac{c F^*}{2 \Delta ^2} \Big] \\
&& +\xi  \Big[k\log8(\Delta+kc)^2+\frac{5\eta^3}{18\Delta^5}\Big(19-\frac{5}{\Delta ^2}\Big)\Big] \\
&&+\frac{5}{4}\gamma\eta  k\log\Big(\frac{1+ks}{1-k s}\Big) 
+\frac{\eta\xi^2}{3\Delta}\Big(8+\frac{4}{\Delta^2}-\frac{3}{\Delta^4}\Big)s \\
&& -\frac{\eta^5s}{4\Delta}\Big(\frac{7}{6\Delta^8}-\frac{17}{3\Delta^6}+\frac{33}{10\Delta^4}+\frac{22}{5\Delta^2}+\frac{44}{5}\Big),
\end{eqnarray*}

\begin{eqnarray*}
\frac{\delta{q}_3}{b} &=& \frac{1}{2k}\frac{c}{\Delta} 
\\
\frac{\delta{q}_4}{b} &=& -\frac{\eta}{6k}\Big(\frac{1}{\Delta ^3}-8 \tilde{E}\Big) 
\\
\frac{\delta{q}_5}{b} &=& \frac{1}{k}\Big\{
\frac{1}{3}(4E^*-F^*)\xi +\frac{c}{4\Delta}\times\\
&& \times\Big[\Big(\frac{1}{\Delta^4}-\frac{1}{3\Delta^2}-\frac{8}{3}\Big)\eta^2
+\Big(\frac{1}{\Delta^2}+4\Big)\xi s\Big] \Big\} 
\\
\frac{\delta{q}_6}{b} &=& \frac{1}{k}\frac{1}{4}\Big\{ \frac{\eta^3}{9\Delta^5}\Big(\frac{5}{\Delta^2}-19\Big)
-k\log8(\Delta+kc)^2 \\
&& +\frac{\eta\xi}{k^2\Delta}\Big(\frac{k^2}{\Delta^4}-\frac{1}{\Delta^2}-2\Big)s \Big\} 
\end{eqnarray*}

\begin{eqnarray*}
\frac{\delta\Phi_1}{\Phi} &=& \frac{1}{\Delta}-2\tilde{K}
\\
\frac{\delta\Phi_2}{\Phi} &=& -\frac{\eta{c}}{\Delta^3}
\\
\frac{\delta\Phi_3}{\Phi} &=& \frac{4}{3}\eta^2(\tilde{E}-\tilde{K})
+\frac{1}{2\Delta^3}\Big[ \Big(3-\frac{1}{\Delta ^2}\Big) \eta^2-\xi  s \Big]
\\
\frac{\delta\Phi_4}{\Phi} &=& \frac{\eta{c}}{2\Delta^5}\Big[
\frac{1}{3}\Big(\frac{5}{\Delta ^2}-11\Big) \eta^2+3 \xi  s\Big] 
\\
\frac{\delta\Phi_5}{\Phi} &=& \gamma\Big[
1-\frac{1}{2\Delta^2}+\tilde{K}\Big(\frac{1}{\Delta }-2\tilde{K}\Big)-\frac{k^2}{2\Delta^3}F^*sc\Big] \\
&&+\xi^2\Big[\frac{1}{\Delta^3}\Big(\frac{1}{2\Delta^2}-1\Big)+\frac{\tilde{K}-4\tilde{E}}{k^2}\Big] \\
&& +\frac{\eta^4}{9}\Big[\frac{1}{8\Delta^5}\Big(\frac{35}{\Delta ^4}-\frac{190}{\Delta^2}+227\Big)
+14\tilde{E}-11\tilde{K}\Big] \\
&& +\frac{\eta^2\xi}{4\Delta^5}\Big(\frac{5}{\Delta^2}-17\Big)s
\\
\frac{\delta\Phi_6}{\Phi} &=& 
\frac{k^2\xi}{\Delta}\Big[\frac{5\eta^3}{9\Delta^6}\Big(19-\frac{7}{\Delta ^2}\Big)c-2\Big]s 
+\frac{c\eta}{\Delta^5} \\
&&\times\Big[\Big(4-\frac{5}{2\Delta^2}\Big)\xi^2 
-\frac{\eta^4}{8\Delta^2}\Big(\frac{7}{\Delta^4}-\frac{98}{3\Delta^2}+\frac{101}{3}\Big)\Big] \\
&&+\frac{\gamma\eta}{2\Delta^2} 
\Big[\Big(2-\frac{4\tilde{K}}{\Delta}+\frac{1}{\Delta^2}\Big)c 
-\frac{F^*}{\Delta}\Big(\frac{k^2}{\Delta ^2}-2\Big)s\Big] 
\end{eqnarray*}
\begin{eqnarray*}
\frac{\delta{Q}_2}{B} &=& -\frac{k{s}}{2\Delta} \\
\frac{\delta{Q}_3}{B} &=& \frac{\eta}{2k}\Big[
\frac{k^2}{\Delta}\Big(\frac{1}{\Delta^2}+1\Big)cs+E^*-F^*\Big]
\\
\frac{\delta{Q}_4}{B} &=& \frac{k}{3}\Big[\Big(\frac{1}{\Delta^3}-8\tilde{E}\Big)\xi
+\frac{1}{4}\frac{s}{\Delta}\Big(\frac{3}{\Delta^4}-\frac{7}{\Delta^2}-14\Big) \eta^2\Big]
\\
\frac{\delta{Q}_5}{B} &=& \frac{k}{9}\Big\{ (14E^*-11F^*)\eta^2
+\frac{3}{\Delta}\Big[ \xi\Big(8+\frac{1}{\Delta^2}-\frac{3}{\Delta^4}\Big) \\
&& +\eta^2\Big(\frac{7}{2}+\frac{11}{4\Delta^2}+\frac{2}{\Delta^4}-\frac{5}{4\Delta^6}\Big)s\Big]c
\Big\}\eta
\\
\frac{\delta{Q}_6}{B} &=& \frac{k}{9}\Big\{ 
\frac{\gamma}{2\Delta}\Big[s\Big(\frac{k^2}{\Delta }-2 \tilde{K}\Big)+\frac{c
   F^*}{2 \Delta ^2}\Big] \\
&&+\frac{\eta^2\xi}{6\Delta^5}\Big(19-\frac{5}{\Delta^2}\Big)
+\frac{\xi^2s}{3\Delta}\Big(2-\frac{k^2}{\Delta ^4}+\frac{1}{\Delta^2}\Big) \\
&&+\frac{1}{4} \gamma  k \log\Big(\frac{1+ks}{1-ks}\Big) \\
&&-\frac{\eta^4 s}{2\Delta}\Big(\frac{35}{72 \Delta
   ^8}-\frac{85}{36 \Delta ^6}+\frac{11}{8 \Delta ^4}+\frac{11}{6 \Delta
   ^2}+\frac{11}{3}\Big)
\Big\}
\end{eqnarray*}

The inverse transformation, needed for the computation of initial conditions in prime variables, is computed from the corrections
\[
\beta'=\beta+\gamma(\Phi)\sum_{1\le{i}\le6}\delta\beta'_i(\phi,q,\Phi,Q),
\]
where the terms $\delta\beta'_i$ are formally the opposite of $\delta\beta_i$, except for
\begin{eqnarray*}
\delta\phi'_5 &=& \gamma\Big(\frac{1}{2\Delta}-\tilde{K}\Big)F^*-\delta\phi_5, 
\\
\delta\phi'_6 &=& \frac{\gamma\eta}{36\Delta}\Big[
\frac{c}{\Delta^2}F^*+74\Big(\frac{1}{\Delta}-2\tilde{K}\Big)s\Big]-\delta\phi_6, 
\\
\frac{\delta\Phi'_5}{\Phi} &=& 2\gamma\tilde{K}\Big(\frac{1}{\Delta}-\tilde{K}\Big)
-\frac{\gamma}{2\Delta^2}\Big(1+\frac{k^2}{\Delta}F^*sc\Big)
-\frac{\delta\Phi_5}{\Phi}
\\
\frac{\delta\Phi'_6}{\Phi} &=& \Big[
\Big(\frac{15}{\Delta^2}+14-\frac{58}{\Delta}\tilde{K}\Big)c+\frac{11}{\Delta}\Big(2-\frac{k^2}{\Delta^2}\Big)F^*s\Big] \\
&& \times\frac{\gamma\eta}{18\Delta^2}-\frac{\delta\Phi_6}{\Phi}
\\
\frac{\delta{Q'}_6}{B} &=& \frac{11k\gamma}{18\Delta}\Big[
\frac{c}{2\Delta^2}F^*+\Big(\frac{1}{\Delta}-2\tilde{K}\Big)s\Big] -\frac{\delta{Q}_6}{B}
\end{eqnarray*}

\color{black}

\section{Lindstedt series solution} \label{a:tables}

Specifically, the Lindstedt series we obtained are conveniently arranged in the form of the following summations:
\begin{eqnarray} \nonumber
n &=& \sum_{m=0}^2\sum_{j=0}^m\sum_{k=0}^{m-j}\alpha^{m-j-k}
 \\ \label{nsL}
&& \times \left(\frac{Q'_0/\Omega}{b}\right)^{2j}\left(\frac{q'_0}{b}\right)^{2k} n_{m,j,k} 
\\ \nonumber 
q' &=& \sum_{m=0}^2\sum_{i=0}^m\sum_{j=0}^m\sum_{k=0}^{m-j}\alpha^{m-j-k}\left(\frac{Q'_0/\Omega}{b}\right)^{2j} \\ \nonumber
&& \times \left(\frac{q'_0}{b}\right)^{2k}\Big[ c_{m,i,j,k}q'_0\cos(2i+1)\Omega\tau \\ \label{qpL}
&& +s_{m,i,j,k}(Q'_0/\Omega)\sin(2i+1)\Omega\tau \Big]
\\ \nonumber 
Q' &=& \sum_{m=0}^2\sum_{i=0}^m\sum_{j=0}^m\sum_{k=0}^{m-j}\alpha^{m-j-k}
\left(\frac{Q'_0/\Omega}{b}\right)^{2j} \\ \nonumber
&& \times\left(\frac{q'_0}{b}\right)^{2k}\Big[ C_{m,i,j,k}Q'_0\cos(2i+1)\Omega\tau \\ \label{QpL}
&& +S_{m,i,j,k}(q'_0\Omega)\sin(2i+1)\Omega\tau \Big]
\\ \nonumber
p &=& \frac{(64/13)k}{(\mathrm{K}-\mathrm{E})^3}\alpha^2\left(\frac{q'_0}{b}\sin\Omega\tau+\frac{Q'_0/\Omega}{b}\cos\Omega\tau\right) \\ \nonumber
&&+ \frac{q'_0}{b}\frac{Q'_0/\Omega}{b}\sum_{m=0}^2\sum_{i=1}^{m+1}\sum_{j=0}^m\sum_{k=0}^{m-j}
\alpha^{m-j-k} \\ \nonumber
&& \times \left(\frac{Q'_0/\Omega}{b}\right)^{2j}\left(\frac{q'_0}{b}\right)^{2k}
\kappa_{m,i,j,k}\cos2i\Omega\tau \\ \nonumber
&&  +\sum_{m=0}^2\sum_{i=1}^{m+1}\sum_{j=0}^{2m+1}\sum_{k=0}^{m+1-j}
\alpha^{m-j-k+1} \\
&& \times\left(\frac{Q'_0/\Omega}{b}\right)^{2j}\left(\frac{q'_0}{b}\right)^{2k}
\sigma_{m,i,j,k}\sin2i\Omega\tau
\\ \nonumber
d &=& \frac{\mathrm{K}}{\mathrm{K}-\mathrm{E}}+
\sum_{m=0}^2\sum_{j=0}^{m+1}\sum_{k=0}^{m+1-j} \alpha^{m-j-k+1} \\ \label{dsL}
&& \times\left(\frac{Q'_0/\Omega}{b}\right)^{2j}\left(\frac{q'_0}{b}\right)^{2k}
d_{m,j,k} 
\end{eqnarray}
in which $\alpha=(\Omega/\omega)^2$, and the exact numeric coefficients $n_{m,j,k}$, $c_{m,i,j,k}$, $s_{m,i,j,k}$, $C_{m,i,j,k}$, $S_{m,i,j,k}$, $\kappa_{m,i,j,k}$, $\sigma_{m,i,j,k}$, and $d_{m,j,k}$, are given in following tables, where the abbreviations $E\equiv{E}(k^2)$, $F\equiv{F}(k^2)$ are used for convenience.


\begin{table}[htb]
\caption{Coefficients $n_{m,j,k}$}
\centering
\begin{tabular}{@{}l@{}}\hline
$n_{0,0,0}=1$ \\ 
$n_{1,0,0}=-\frac{-64{E}{K}-256 {K}^2+320{E}^2+63}{144 ({E}-{K})^2}$ \\
$n_{1,0,1}=n_{1,1,0}=\frac{3(11 {K}-14 {E})}{64({K}-{E})}$ \\
$n_{2,0,0}=-\frac{2(4 {E}-{K})^2}{81({E}-{K})^2}$ \\ 
$n_{2,0,1}=-\frac{-370 {E}{K}+35 {K}^2+344 {E}^2}{96({E}-{K})^2}$ \\ 
$n_{2,0,2}=\frac{-12892 {E}{K}+5459 {K}^2+7244 {E}^2}{16384({E}-{K})^2}$ \\ 
$n_{2,1,0}=-\frac{-162 {E}{K}+19 {K}^2+152 {E}^2}{32({E}-{K})^2}$ \\ 
$n_{2,1,1}=2n_{2,2,0}=-\frac{-1892 {E}{K}+349 {K}^2+2164 {E}^2}{8192({E}-{K})^2}$ \\ 
\hline
\end{tabular}
\label{t:nmjk}
\end{table}

\begin{table}[htb]
\caption{Coefficients $c_{m,i,j,k}$}
\centering
\begin{tabular}{@{}l@{}}\hline
$c_{0,0,0,0}=1$, \qquad 
$c_{1,0,0,0}=c_{1,1,0,0}=0$ \\ 
$c_{2,0,0,0}=c_{2,1,0,0}=c_{2,2,0,0}=c_{2,2,0,1}=c_{2,2,1,0}=0$ \\ 
$c_{1,0,1,0}=-c_{1,1,1,0}=\frac{3}{4}$, \qquad 
$c_{1,1,0,1}=-c_{1,0,0,1}=\frac{1}{4}$ \\ 
$c_{2,0,0,1}=-c_{2,1,0,1}=\frac{12 {E}-{K}}{16({E}-{K})}$ \\ 
$c_{2,0,0,2}=-\frac{-9724 {E}{K}+4451 {K}^2+4652 {E}^2}{196608({E}-{K})^2}$ \\ 
$c_{2,0,1,0}=-\frac{-526 {E}{K}+47 {K}^2+488 {E}^2}{192({E}-{K})^2}$ \\ 
$c_{2,0,1,1}=-\frac{5 (-7964{E} {K}+2971 {K}^2+5452 {E}^2)}{98304({E}-{K})^2}$ \\ 
$c_{2,0,2,0}=-\frac{-26972 {E}{K}+9019 {K}^2+21004 {E}^2}{196608({E}-{K})^2}$ \\ 
$c_{2,1,0,2}=\frac{-286 {E}{K}+137 {K}^2+122 {E}^2}{8192({E}-{K})^2}$ \\ 
$c_{2,1,1,0}=\frac{-526 {E}{K}+47 {K}^2+488 {E}^2}{192({E}-{K})^2}$ \\ 
$c_{2,1,1,1}=\frac{5 (-1804 {E}{K}+689 {K}^2+1196 {E}^2)}{16384({E}-{K})^2}$ \\ 
$c_{2,1,2,0}=\frac{3 (-352 {E}{K}+89 {K}^2+344 {E}^2)}{16384({E}-{K})^2}$ \\ 
$c_{2,2,0,2}=\frac{-2860 {E} {K}+1163 {K}^2+1724 {E}^2}{196608({E}-{K})^2}$ \\ 
$c_{2,2,1,1}=-\frac{5(-2860 {E} {K}+1163 {K}^2+1724{E}^2)}{98304({E}-{K})^2}$ \\ 
$c_{2,2,2,0}=\frac{5 (-2860 {E}{K}+1163 {K}^2+1724 {E}^2)}{196608({E}-{K})^2}$ \\
   \hline
\end{tabular}
\label{t:cmijk}
\end{table}

\begin{table}[htb]
\caption{Coefficients $s_{m,i,j,k}$}
\centering
\begin{tabular}{@{}l@{}}\hline
$s_{0,0,0,0}=-1$ \\ 
$s_{1,0,0,0}=\frac{128(4{E}-{K})}{9 (14{E}-11{K})}-\frac{32{E}{K}-96{K}^2+64{E}^2+21}{48({E}-{K})^2}$ \\ 
$s_{1,0,0,1}=\frac{21}{4}$, \quad
$s_{1,0,1,0}=\frac{9}{4}$, \quad 
$s_{1,1,0,1}=-\frac{3}{4}$, \quad
$s_{1,1,1,0}=\frac{1}{4}$ \\ 
$s_{1,1,0,0}=s_{2,1,0,0}=s_{2,2,0,0}=s_{2,2,0,1}=s_{2,2,1,0}=0$ \\ 
$s_{2,0,0,0}=-n_{2,0,0}$ \\ 
$s_{2,0,0,1}=\frac{-530 {E}{K}-5 {K}^2+472 {E}^2}{384({E}-{K})^2}$ \\ 
$s_{2,0,0,2}=\frac{11396 {E}{K}+6563 {K}^2-34132 {E}^2}{196608({E}-{K})^2}$ \\ 
$s_{2,0,1,0}=-\frac{-278 {E}{K}+41 {K}^2+264 {E}^2}{128({E}-{K})^2}$ \\ 
$s_{2,0,1,1}=-\frac{7 (-17644{E} {K}+6143 {K}^2+13148 {E}^2)}{98304({E}-{K})^2}$ \\ 
$s_{2,0,2,0}=-\frac{-95524 {E}{K}+34373 {K}^2+68468 {E}^2}{196608({E}-{K})^2}$ \\ 
$s_{2,1,0,1}=\frac{-994{E} {K}+83 {K}^2+920 {E}^2}{384({E}-{K})^2}$ \\ 
$s_{2,1,0,2}=\frac{3 (-968 {E}{K}+331 {K}^2+736 {E}^2)}{16384({E}-{K})^2}$ \\ 
$s_{2,1,1,0}=\frac{1}{4}n_{2,0,1}$ \\ 
$s_{2,1,1,1}=-\frac{-5764 {E}{K}+2363 {K}^2+3428 {E}^2}{16384({E}-{K})^2}$ \\ 
$s_{2,1,2,0}=-\frac{-319 {E}{K}+113 {K}^2+233 {E}^2}{4096({E}-{K})^2}$ \\ 
$s_{2,2,0,2}=-5c_{2,2,0,2}$ \\ 
$s_{2,2,1,1}=10c_{2,2,0,2}$, \quad 
$s_{2,2,2,0}=-c_{2,2,0,2}$ \\
\hline
\end{tabular}
\label{t:smijk}
\end{table}

\begin{table}[htb]
\caption{Coefficients $C_{m,i,j,k}$}
\centering
\begin{tabular}{@{}l@{}}\hline
$C_{0,0,0,0}=1$, \quad 
$C_{1,0,0,0}=C_{1,1,0,0}=C_{2,0,0,0}=0$ \\ 
$C_{2,1,0,0}=C_{2,2,0,0}=C_{2,2,0,1}=C_{2,2,1,0}=0$ \\ 
$C_{1,1,0,1}=-C_{1,0,0,1}=\frac{9}{4}$, \quad
$C_{1,0,1,0}=-C_{1,1,1,0}=\frac{3}{4}$ \\ 
$C_{2,0,0,1}=-C_{2,1,0,1}=-3c_{2,0,0,1}$ \\ 
$C_{2,0,0,2}=-\frac{-66748 {E}{K}+32531 {K}^2+27116 {E}^2}{196608({E}-{K})^2}$ \\ 
$C_{2,0,1,0}=-C_{2,1,1,0}=\frac{-98 {E}{K}+{K}^2+88 {E}^2}{192({E}-{K})^2}$ \\ 
$C_{2,0,1,1}=-\frac{-65516 {E}{K}+26527 {K}^2+39772 {E}^2}{98304({E}-{K})^2}$ \\ 
$C_{2,0,2,0}=c_{2,0,2,0}$ \\ 
$C_{2,1,0,2}=\frac{9 (11 {E}{K}+8 {K}^2-37 {E}^2)}{4096({E}-{K})^2}$ \\ 
$C_{2,1,1,1}=\frac{3 (-7612 {E}{K}+3089 {K}^2+4604 {E}^2)}{16384({E}-{K})^2}$ \\ 
$C_{2,1,2,0}=c_{2,1,2,0}$, \qquad
$C_{2,2,0,2}=25c_{2,2,0,2}$ \\ 
$C_{2,2,1,1}=-50c_{2,2,0,2}$, \qquad
$C_{2,2,2,0}=5c_{2,2,0,2}$ \\ 
\hline
\end{tabular}
\label{t:Cmijk}
\end{table}

\begin{table}[htb]
\caption{Coefficients $S_{m,i,j,k}$}
\centering
\begin{tabular}{@{}l@{}}\hline
$S_{0,0,0,0}=1$, \qquad
$S_{1,0,0,0}=s_{1,0,0,0}$ \\ 
$S_{1,0,0,1}=\frac{11}{15}S_{1,0,1,0}=\frac{11}{4}$, \quad
$S_{1,1,0,1}=-\frac{1}{3}S_{1,1,1,0}=\frac{3}{4}$ \\
$S_{1,1,0,0}=S_{2,1,0,0}=S_{2,2,0,0}=S_{2,2,0,1}=S_{2,2,1,0}=0$ \\ 
$S_{2,0,0,0}=3s_{2,0,0,0}$ \\ 
$S_{2,0,0,1}=\frac{-1574 {E}{K}+193 {K}^2+1480 {E}^2}{1152({E}-{K})^2}$ \\ 
$S_{2,0,0,2}=\frac{-133892 {E}{K}+56701 {K}^2+75220 {E}^2}{196608({E}-{K})^2}$ \\ 
$S_{2,0,1,0}=-\frac{-1226 {E}{K}+127 {K}^2+1144 {E}^2}{384({E}-{K})^2}$ \\ 
$S_{2,0,1,1}=-\frac{-51436 {E}{K}+14687 {K}^2+46172 {E}^2}{98304({E}-{K})^2}$ \\ 
$S_{2,0,2,0}=-\frac{-16412 {E}{K}+139 {K}^2+25804 {E}^2}{196608({E}-{K})^2}$ \\ 
$S_{2,1,0,1}=-\frac{1}{4}n_{2,0,1}$ \\ 
$S_{2,1,0,2}=\frac{3 (-1496 {E}{K}+637 {K}^2+832 {E}^2)}{16384({E}-{K})^2}$ \\ 
$S_{2,1,1,0}=-\frac{-254 {E}{K}+13 {K}^2+232 {E}^2}{128({E}-{K})^2}$ \\ 
$S_{2,1,1,1}=\frac{3 (-7172 {E}{K}+2719 {K}^2+4804 {E}^2)}{16384({E}-{K})^2}$ \\ 
$S_{2,1,2,0}=-9c_{2,1,0,2}$ \\ 
$S_{2,2,0,2}=\frac{1}{5}S_{2,2,2,0}=-\frac{1}{10}S_{2,2,1,1}=C_{2,2,2,0}$ \\ 
\hline
\end{tabular}
\label{t:Smijk}
\end{table}

\begin{table}[htb]
\caption{Coefficients $\kappa_{m,i,j,k}$}
\centering
\begin{tabular}{@{}l@{}}\hline
$\kappa_{0,1,0,0}=\frac{3}{4}$ \\ 
$\kappa_{1,1,0,0}=-\frac{7 (64 {E}{K}-128 {K}^2+64 {E}^2+27)}{192({E}-{K})^2}$ \\ 
$\kappa_{1,1,0,1}=\frac{5}{12}n_{1,0,1}$, \qquad
$\kappa_{1,1,1,0}=\frac{17}{12}n_{1,0,1}$ \\ 
$\kappa_{1,2,0,0}=\kappa_{2,2,0,0}=\kappa_{2,3,0,0}=\kappa_{2,3,0,1}=\kappa_{2,3,1,0}=0$ \\ 
$\kappa_{1,2,0,1}=-\kappa_{1,2,1,0}=\frac{13}{24}n_{1,0,1}$ \\ 
$\kappa_{2,1,0,0}=-\frac{9}{4}n_{2,0,0}$ \\ 
$\kappa_{2,1,0,1}=-\frac{-3446 {E}{K}+337 {K}^2+3208 {E}^2}{1152({E}-{K})^2}$ \\ 
$\kappa_{2,1,0,2}=\frac{7 (-2948 {E}{K}+1237 {K}^2+1684 {E}^2)}{262144({E}-{K})^2}$ \\ 
$\kappa_{2,1,1,0}=-3S_{2,0,0,1}$ \\ 
$\kappa_{2,1,1,1}=-\frac{-153604 {E}{K}+52853 {K}^2+115988 {E}^2}{131072({E}-{K})^2}$ \\ 
$\kappa_{2,1,2,0}=\frac{9 (-4796 {E}{K}+3067 {K}^2+172 {E}^2)}{262144({E}-{K})^2}$ \\ 
$\kappa_{2,2,0,1}=\frac{-6842{E} {K}+559 {K}^2+6328 {E}^2}{2304({E}-{K})^2}$ \\ 
$\kappa_{2,2,0,2}=-\kappa_{2,2,2,0}=\frac{27 (-572 {E}{K}+251 {K}^2+300 {E}^2)}{32768({E}-{K})^2}$ \\ 
$\kappa_{2,2,1,0}=\frac{7}{6}S_{2,1,1,0}$, \qquad
$\kappa_{2,2,1,1}=\frac{416}{25}\kappa_{1,1,0,1}^2$ \\ 
$\kappa_{2,3,0,2}=\kappa_{2,3,2,0}=\frac{-202268 {E} {K}+81907 {K}^2+122764 {E}^2}{786432({E}-{K})^2}$ \\ 
$\kappa_{2,3,1,1}=-\frac{10}{3}\kappa_{2,3,0,2}$ \\ 
\hline
\end{tabular}
\label{t:kamijk}
\end{table}

\begin{table}[htb]
\caption{Coefficients $\sigma_{m,i,j,k}$}
\centering
\begin{tabular}{@{}l@{}}\hline
$\sigma_{0,1,0,0}=\sigma_{1,1,0,0}=\sigma_{1,2,0,0}=\sigma_{1,2,0,1}=\sigma_{1,2,1,0}=0$ \\ 
$\sigma_{2,1,0,0}=\sigma_{2,1,1,0}=\sigma_{2,1,1,1}=\sigma_{2,2,0,0}=\sigma_{2,2,0,1}=0$ \\ 
$\sigma_{2,2,1,0}=\sigma_{2,3,0,0}=\sigma_{2,3,0,1}=\sigma_{2,3,0,2}=\sigma_{2,3,1,0}=0$ \\ 
$\sigma_{2,3,1,1}=\sigma_{2,3,2,0}=0$, \qquad
$\sigma_{0,1,0,1}=-\sigma_{0,1,1,0}=\frac{3}{8}$ \\ 
$\sigma_{1,1,0,1}=-\frac{128 {E}^2+176 {K} {E}-304{K}^2+63}{96 ({E}-{K})^2}$ \\ 
$\sigma_{1,1,0,2}=\frac{5}{6}n_{1,0,1}$ \\ 
$\sigma_{1,1,1,0}=\frac{64 {E}^2+32{K} {E}-96 {K}^2+21}{64({E}-{K})^2}$ \\ 
$\sigma_{1,1,1,1}=\frac{3}{2}n_{1,0,1}$, \qquad 
$\sigma_{1,1,2,0}=-\frac{1}{3}n_{1,0,1}$ \\ 
$\sigma_{1,2,0,2}=\frac{1}{4}\kappa_{1,2,0,1}$, \qquad
$\sigma_{1,2,1,1}=-6\sigma_{1,2,0,2}$ \\ 
$\sigma_{1,2,2,0}=\sigma_{1,2,0,2}$, \qquad
$\sigma_{2,1,0,1}=-3n_{2,0,0}$ \\ 
$\sigma_{2,1,0,2}=\frac{8(4{E}-{K})}{15({E}-{K})}\kappa_{1,1,0,1}$ \\ 
$\sigma_{2,1,0,3}=\frac{298364{E}^2-510796 {K} {E}+211631 {K}^2}{524288({E}-{K})^2}$ \\ 
$\sigma_{2,1,1,2}=\frac{136636 {E}^2-358028 {K} {E}+178111 {K}^2}{524288({E}-{K})^2}$ \\ 
$\sigma_{2,1,2,0}=\frac{4}{3}\sigma_{2,1,0,2}$ \\ 
$\sigma_{2,1,2,1}=\frac{3 (138092{E}^2-233596 {K} {E}+96107 {K}^2)}{524288({E}-{K})^2}$ \\ 
$\sigma_{2,1,3,0}=\frac{3 (15916{E}^2-17468 {K} {E}+4891 {K}^2)}{524288({E}-{K})^2}$ \\ 
$\sigma_{2,2,0,2}=\frac{392 {E}^2-422 {K} {E}+39 {K}^2}{512({E}-{K})^2}$ \\ 
$\sigma_{2,2,0,3}=\frac{10420 {E}^2-17732{K} {E}+7321 {K}^2}{65536({E}-{K})^2}$ \\ 
$\sigma_{2,2,1,1}=-\frac{700{E}^2-761 {K} {E}+52 {K}^2}{192({E}-{K})^2}$ \\ 
$\sigma_{2,2,1,2}=\frac{3 (9812{E}^2-13156 {K} {E}+4577 {K}^2)}{65536({E}-{K})^2}$ \\ 
$\sigma_{2,2,2,0}=\frac{2072 {E}^2-2290{K} {E}+65 {K}^2}{4608({E}-{K})^2}$ \\ 
$\sigma_{2,2,2,1}=-\frac{3 (20764{E}^2-34892 {K} {E}+14299 {K}^2)}{65536({E}-{K})^2}$ \\ 
$\sigma_{2,2,3,0}=\frac{3 (76{E}^2-572 {K} {E}+343 {K}^2)}{65536({E}-{K})^2}$ \\ 
$\sigma_{2,3,0,3}=-\sigma_{2,3,3,0}=\frac{1}{6}\kappa_{2,3,0,2}$ \\ 
$\sigma_{2,3,1,2}=-\sigma_{2,3,2,1}=-\frac{5}{2}\kappa_{2,3,0,2}$ \\ 
\hline
\end{tabular}
\label{t:simijk}
\end{table}

\begin{table}[htb]
\caption{Coefficients $d_{m,j,k}$}
\centering
\begin{tabular}{@{}l@{}}\hline
$d_{0,0,0}=-\frac{1-4{K}^2}{({E}-{K})^2}$ \\ 
$d_{0,0,1}=d_{0,1,0}=\frac{3}{4}$, \qquad
$d_{1,0,0}=d_{2,0,0}=d_{2,1,0}=0$ \\ 
$d_{1,0,1}=-\frac{16 {E} {K}-272{K}^2+256 {E}^2+63}{48({E}-{K})^2}$ \\ 
$d_{1,0,2}=\frac{9}{8}n_{1,0,1}$ \\ 
$d_{1,1,0}=-\frac{7(-32 {E} {K}-32 {K}^2+64 {E}^2+9)}{96({E}-{K})^2}$ \\ 
$d_{1,1,1}=2d_{1,2,0}=\frac{1}{4}n_{1,0,1}$, \qquad 
$d_{2,0,1}=6n_{2,0,0}$ \\ 
$d_{2,0,2}=-\frac{-1730 {E}{K}+205 {K}^2+1624 {E}^2}{384({K}-{E})^2}$ \\ 
$d_{2,0,3}=\frac{-55220 {E}{K}+22837 {K}^2+32356 {E}^2}{49152({K}-{E})^2}$ \\ 
$d_{2,1,1}=-\frac{-799 {E}{K}+68 {K}^2+740 {E}^2}{48({K}-{E})^2}$ \\ 
$d_{2,1,2}=\frac{-1364 {E}{K}+733 {K}^2+388 {E}^2}{16384({K}-{E})^2}$ \\ 
$d_{2,2,0}=-\frac{-1382 {E}{K}+139 {K}^2+1288 {E}^2}{384({K}-{E})^2}$ \\ 
$d_{2,2,1}=3d_{2,3,0}=-\frac{-2332 {E}{K}+719 {K}^2+1964 {E}^2}{16384({K}-{E})^2}$ \\ 
\hline
\end{tabular}
\label{t:dmjk}
\end{table}

\end{document}